\documentclass[11pt]{amsart}
\usepackage{amsmath}
\usepackage{amsxtra}
\usepackage{amscd}
\usepackage{amsthm}
\usepackage{amsfonts}
\usepackage{amssymb}
\usepackage{eucal}

\numberwithin{equation}{subsection}

\begin{document}

\newcommand{\thmref}[1]{Theorem~\ref{#1}}
\newcommand{\secref}[1]{Section~\ref{#1}}
\newcommand{\defref}[1]{Definition~\ref{#1}}
\newcommand{\lemref}[1]{Lemma~\ref{#1}}
\newcommand{\propref}[1]{Proposition~\ref{#1}}
\newcommand{\corref}[1]{Corollary~\ref{#1}}
\newcommand{\remref}[1]{Remark~\ref{#1}}
\newcommand{\nc}{\newcommand}

\nc{\lrarrow}{\longrightarrow}
\nc{\mc}{\mathcal}
\nc{\on}{\operatorname}
\nc{\Wick}{{\mb :}}
\nc{\delz}{\partial_z}
\nc{\ddz}{\frac{\partial}{\partial z}}
\nc{\Z}{{\mbb Z}}
\nc{\C}{{\mbb C}}
\nc{\Oo}{{\mc O}}
\nc{\D}{{\mc D}}
\nc{\E}{{\mc E}}
\nc{\bib}{\bibitem}
\nc{\pone}{\Pro^1}
\nc{\pa}{\partial}
\nc{\K}{{\mc K}}
\nc{\W}{{\mc W}}
\nc{\Ll}{{\mc L}}
\nc{\la}{\lambda}
\nc{\ep}{\epsilon}
\nc{\su}{\widehat{{\mf s}{\mf l}}_2}
\nc{\sw}{{\mf s}{\mf l}}
\nc{\g}{{\mf g}}
\nc{\h}{{\mf h}}
\nc{\n}{{\mf n}}
\nc{\ab}{\mf{a}}
\nc{\is}{{\mb i}}
\nc{\V}{\mc{V}}
\nc{\Vir}{\mc{V}ir}
\nc{\inv}{^{-1}}
\nc{\vac}{|0\rangle}
\nc{\ol}{\overline}
\nc{\wt}{\widetilde}
\nc{\wh}{\widehat}
\nc{\dst}{\displaystyle}
\nc{\al}{\alpha}

\nc{\mb}{\mathbf}
\nc{\mf}{\mathfrak}
\nc{\ot}{\otimes}
\nc{\mbb}{\mathbb}
\nc{\ghat}{\wh{\g}}
\nc{\gtil}{\wt{\g}}
\nc{\Gtil}{\wt{G}}
\newcommand\Cx{{\mbb C}^\times}
\nc{\un}{\underline}
\nc{\BB}{{\mc B}}
\nc{\bb}{{\mf b}}
\nc{\kk}{{\mf k}}
\nc{\khat}{\wh{\kk}}
\nc{\Y}{{\mc Y}}
\nc{\frob}{\times}
\nc{\sm}{\setminus}
\nc{\bs}{\backslash}
\nc{\Pone}{{\mathbb P}^1}
\nc{\Aone}{{\mathbb A}^1}
\nc{\AutO}{\on{Aut}\Oo}
\nc{\Der}{{\on{Der}\,}}
\nc{\DerO}{{\on{Der}\Oo}}
\nc{\AutK}{\on{Aut}\K}
\nc{\Aut}{\on{Aut}}
\nc{\Sug}{{\mbb S}}
\nc{\Mg}{{{\mf M}_g}}
\nc{\hv}{h^{\vee}}
\nc{\Proj}{{\mc P}roj\,}
\nc{\Conn}{{\mc C}onn\,}
\nc{\ExConn}{{\mc E}x{\mc C}onn}
\nc{\Op}{{\mc Op}}
\nc{\Hitch}{\on{Hitch}}
\nc{\Higgs}{\on{Higgs}}
\nc{\Gm}{{\mathbb G}_m}
\nc{\Xg}{{\mf X}_g}
\nc{\Xf}{{\mf X}}
\nc{\AutXg}{{\mc Aut}_{{\mf X}_g}}
\nc{\T}{{\mc T}}
\nc{\boxt}{\boxtimes}
\nc{\Mgo}{{\mf M}_{g,1}}
\nc{\x}{{\mf x}}
\nc{\Pc}{{\mc P}}
\nc{\Pf}{{\mf P}}
\nc{\Ac}{{\mc A}}
\nc{\sug}{{\mf s}}
\nc{\Gr}{{\mc Gr}}
\nc{\Hh}{{\mc H}}
\nc{\Vgk}{V_{\g,K}}
\nc{\Vck}{Vir_{c_k}^{2n}}
\nc{\Sym}{\on{Sym}}
\nc{\M}{{\mf M}}
\nc{\Mcal}{{\mc M}}
\nc{\Mhat}{\wh{\M}}
\nc{\Mghat}{\wh{\M}_{g,1}}
\nc{\Mtil}{\wt{\M}}
\nc{\Cc}{{\mc C}}
\nc{\vphi}{\varphi}
\nc{\unVir}{Vir^{0}}
\nc{\unV}{\un{V}}
\nc{\UnVir}{{\mc V}ir^{0}}
\nc{\UnV}{\un{\mc V}}
\nc{\unsug}{\un{\sug}}
\nc{\unSug}{\un{\Sug}}
\nc{\mm}{{\mf m}}
\nc{\Gn}{G_n(\Oo)}
\nc{\An}{\Aut_n\Oo}
\nc{\Atn}{\Aut_{2n}\Oo}

\nc{\Bun}{\on{Bun}}
\nc{\und}{\underset}

\title{Geometric Realization of the Segal--Sugawara Construction}

\author{David Ben-Zvi}

\address{Department of Mathematics, University of Chicago, Chicago IL
60637}

\author{Edward Frenkel}\thanks{Partially supported by grants from the
Packard Foundation and the NSF}

\address{Department of Mathematics, University of California,
Berkeley, CA 94720}

\dedicatory{To Graeme Segal on his 60th birthday}

\date{January 2003}

\begin{abstract}
We apply the technique of localization for vertex algebras to the
Segal-Sugawara construction of an ``internal'' action of the Virasoro
algebra on affine Kac-Moody algebras. The result is a lifting of
twisted differential operators from the moduli of curves to the moduli
of curves with bundles, with arbitrary decorations and complex
twistings. This construction gives a uniform approach to
a collection of phenomena describing the geometry of the moduli spaces
of bundles over varying curves: the KZB equations and heat kernels on
non-abelian theta functions, their critical level limit giving the
quadratic parts of the Beilinson-Drinfeld quantization of the Hitchin
system, and their infinite level limit giving a Hamiltonian
description of the isomonodromy equations.
\end{abstract}

\maketitle

\section{Introduction.}

\subsection{Uniformization.}

Let $G$ be a complex connected simply-connected simple algebraic group
with Lie algebra $\g$, and $X$ a smooth projective curve over
$\C$. The geometry of $G$--bundles on $X$ is intimately linked to
representation theory of the affine Kac-Moody algebra $\ghat$, the
universal central extension of the loop algebra $L\g=\g\ot\K$, where
$\K=\C((t))$. More precisely, let $G(\Oo)$, where $\Oo=\C[[t]]$,
denote the positive half of the loop group $G(\K)$. There is a
principal $G(\Oo)$--bundle over the moduli space (stack) $\Bun_G(X)$
of $G$--bundles on $X$, which carries an infinitesimally simply
transitive action of $L\g$. This provides an infinitesimal
``uniformization'' of the moduli of $G$--bundles. Moreover, this
uniformization lifts to an infinitesimal action of $\ghat$ on the
``determinant'' line bundle $\Cc$ on $\Bun_G(X)$, whose sections give
the nonabelian versions of the spaces of theta functions.

The geometry of the moduli space $\Mgo$ of smooth pointed curves of
genus $g$ is similarly linked to the Virasoro algebra $Vir$, the
universal central extension of the Lie algebra $\Der \K$ of
derivations of $\K$ (or vector fields on the punctured disc). Let
$\AutO$ denote the group of automorphisms of the disc.  Then there is
a principal $\AutO$--bundle over $\Mgo$ which carries an
infinitesimally transitive action of $\Der \K$. This uniformization
lifts to an infinitesimal action of the Virasoro algebra on the Hodge
line bundle $\Hh$ on $\Mgo$.

These uniformizations can be used, following \cite{Jantzen,BS}, to
construct {\em localization functors} from representations of the
Virasoro and Kac-Moody algebras to sheaves on the corresponding moduli
spaces. Using these localization functors, one can describe sheaves of
modules over (twisted) differential operators, as well as sheaves on
the corresponding (twisted) cotangent bundles over the moduli
spaces. The sheaves of twisted differential operators and twisted
symbols (functions on twisted cotangent bundles) themselves are
particularly easy to describe in this fashion: they are the
localizations of the vacuum modules of the respective Lie algebras.

The Virasoro algebra acts on the Kac-Moody algebra via its natural
action by derivations on the space $\K$ of Laurent series. Therefore
we form the semidirect product $\gtil = Vir \ltimes \ghat$, which
uniformizes the moduli stack $\Bun_{G,g,1}$ classifying $G$--bundles
over varying pointed curves of genus $g$. There is now a two-parameter
family of line bundles $\Ll_{k,c}$ on which the (two-parameter)
central extension $\gtil$ acts, and the corresponding two-parameter
families of sheaves of twisted differential operators and twisted
symbols. Thus it is natural to look to the structure of $\gtil$ for
information about the variation of $G$--bundles over the moduli of
curves.

\subsection{The Segal-Sugawara construction.}

The interaction between the variation of curves and that of bundles on
curves is captured by a remarkable feature of affine Kac-Moody
algebras, best seen from the perspective of the theory of vertex
algebras. The Segal-Sugawara construction presents the action of the
Virasoro algebra on $\ghat$ as an {\em internal} action, by
identifying infinite quadratic expressions in the Kac-Moody generators
which satisfy the Virasoro relations. More precisely, the Virasoro
algebra embeds in a certain completion of $U\ghat$, and hence acts on
all {\em smooth} representations of $\ghat$ of non-critical
levels. From the point of view of vertex algebras, the construction
simply involves the identification of a conformal structure, i.e., a
certain distinguished vector, in the vacuum module over $\ghat$.

Our objective in this paper is to draw out the different geometric
consequences of this construction in a simple uniform fashion. Our
approach involves four steps:
\begin{enumerate}
\item[$\bullet$] We consider the Segal-Sugawara construction as
defining a $G(\Oo)$--invariant embedding of the vacuum module of the
Virasoro algebra into the vacuum module of the Virasoro--Kac-Moody
algebra $\gtil$.

\item[$\bullet$] Twists of the vacuum modules form sheaves of algebras
on the relevant moduli spaces, and the above construction gives rise
to a homomorphism of these sheaves of algebras.

\item[$\bullet$] The sheaves of twisted differential operators on the
moduli spaces are quotients of the above algebras, and the vertex
algebraic description of the Segal-Sugawara operators guarantees that
the homomorphism descends to a homomorphism between sheaves of
differential operators.

\item[$\bullet$] All of the constructions vary flatly with respect to
the Kac-Moody level and the Virasoro central charge and possess
various ``classical'' limits which may be described in terms of the
twisted cotangent bundles to the respective moduli spaces.
\end{enumerate}

Our main result is the following theorem.  Let $\Pi:\Bun_{G,g}\to\Mg$
denote the projection from moduli of curves and bundles to moduli of
curves, $k,c\in\C$ (the level and charge), $\mu_\g=\hv\on{dim}\g$ and
$c_k=c-\dfrac{k \on{dim}\g}{k+\hv}$. We denote by $\D_{k,c}$,
$\D_{k/\Mg}$, and $\D_{c_k}$ respectively the corresponding sheaves of
twisted differential operators on $Bun_{G,g}$, relative to $\Mg$, and
on $\Mg$. The sheaves of twisted symbols (functions on twisted
cotangent bundles) with twists $\la,\mu$ on $\M_g$ and $Bun_{G,g}$ are
denoted by $\Oo(T^*_{\la,\mu}Bun_{G,g})$ and $\Oo(T^*_{\la\mu}\Mg)$,
respectively.

\subsubsection{\bf Theorem.}\label{main}
{\em
\begin{enumerate}
\item Let $k,c\in\C$ with $k$ not equal to minus the dual Coxeter
number $\hv$.  There is a canonical homomorphism of sheaves of
algebras on $\Mg$, $$\D_{c_k}\longrightarrow \Pi_*\D_{k,c},$$ and an
isomorphism of twisted $\D$--modules on $\Bun_{G,g}$
$$
\D_{k/\Mg}\und{\Oo}\ot \Pi^*\D_{c_k}\cong\D_{k,c}.$$
\item Let $c\in\C$, $k=-\hv$. There is an algebra homomorphism
  $$\Oo(T^*_{\mu_\g}\M) \longrightarrow \Pi_*\D_{-h^{\vee},c}.$$
\item For every $\la,\mu\in\C$ there is a ($\la$--Poisson)
homomorphism
$$\Oo(T^*_{\la\mu}\Mg)\longrightarrow
\Pi_*\Oo(T^*_{\la,\mu}\Bun_{G,g}).$$
\end{enumerate}
}

\subsubsection{}
Moreover, the homomorphisms (2) and (3) are suitably rescaled limits of
(1). In fact we prove a stronger result, valid for moduli of curves
and bundles with arbitrary ``decorations''. Namely, we consider moduli
of curves with an arbitrary number of marked points, jets of local
coordinates and jets of bundle trivializations at these points.  (For
the sake of simplicity of notation, we work with the case of a single
marked point or none -- the multipoint extension is straightforward.)

\subsection{Applications.}
The homomorphism (1) allows us to lift vector fields on $\Mg$ to
differential operators on $\Bun_{G,g}$, which are first order along
the moduli of curves and second order along moduli of bundles. In the
special case of $k\in\Z_+$ and no decorations, this gives a direct
algebraic construction of the heat operators acting on non-abelian
theta functions (the global sections of $\Ll_{k,c}$ along $\Bun_G(X)$)
constructed by different means in \cite{ADW,Hitchin connection,BK,Fa}.
Our approach via vertex algebras makes clear the compatibility of this
heat equation with the projectively flat connection on the bundle of
conformal blocks coming from conformal field theory
(\cite{TUY,FrS})--- the Knizhnik--Zamolodchikov--Bernard equations, or
higher genus generalization of the Knizhnik--Zamolodchikov equations,
\cite{KZ}. This compatibility has also been established in
\cite{Laszlo}. The important feature of our proof which follows
immediately from part (1) of the above theorem is that the heat
operators depend on two {\em complex} parameters $k$ and $c$, allowing
us to consider various limits. This answers some of the questions
raised by Felder in his study of of the general KZB equations over the
moduli space of pointed curves for arbitrary complex levels (see
\cite{Fel}, \S~6).

The geometric significance of the two classical limit statements (2),
(3) comes from the identification of the twisted cotangent bundles of
$\Mg$ and $\Bun_{G}(X)$ as the moduli spaces $T^*_1\M_g=\Proj_g$ of
algebraic curves with projective structures and
$T^*_1Bun_G(X)=\Conn_G(X)$ of $G$--bundles with holomorphic
connections, respectively.  As the parameter $k$ approaches the {\em
  critical level} $-\hv$, the component of the heat kernels along the
moduli of curves drops out, and we obtain commuting global second
order operators on the moduli $\Bun_G(X)$ (for a fixed curve $X$) from
linear functions on the space of projective structures on $X$. Thus we
recover the quadratic part of the Beilinson--Drinfeld quantization of
the Hitchin hamiltonians (see \cite{Hecke}) in part (2).  When we take
$k,c$ to infinity in different fashions, we obtain the classical
limits in part (3). For $\la,\mu=0$ (i.e., as the symbols of our
operators) we recover the quadratic part of the Hitchin integrable
system \cite{Hitchin system}, namely, the hamiltonians corresponding
to the trace of the square of a Higgs field. For $\la\neq 0$ we obtain
a canonical (non-affine) decomposition of the twisted cotangent bundle
$T^*_{\la,\mu}\Bun_G=\ExConn^{\la,\mu}_{G,g}$, the space of {\em
  extended connections} (\cite{BS,szego}), into a product over $\M_g$
of the spaces of projective structures and flat connections. This
recovers a construction of Bloch--Esnault and Beilinson \cite{BE}.
(This decomposition is also described from the point of view of kernel
functions in \cite{thetas}, where it is related to the Klein
construction of projective connections from theta functions and used
to give a conjectural coordinate system on the moduli of bundles.)

For $\la\neq 0$, the homomorphism (3) gives a Hamiltonian action of
the vector fields on $\Mg$ on the moduli space of flat
$G$--bundles. This gives a time-dependent Hamiltonian system (see
\cite{M}), or flat symplectic connection, on the moduli $\Conn_{G,g}$
of flat connections, with times given by the moduli of decorated
curves. By working with arbitrary level structures, we obtain a
similar statement for moduli spaces of connections with {\em arbitrary}
(regular or irregular) singularities. We show that this is in fact a
Hamiltonian description of the equations of isomonodromic deformation
of meromorphic connections with respect to variation of the decorated
curve. In particular, since all of the maps of the theorem come from a
single multi--parametric construction, \thmref{main} immediately
implies the following picture suggested by the work of \cite{LO} (see
also \cite{Ivanov}), with arrows representing degenerations:
\begin{equation*}
\begin{array}{ccc}
\mbox{Quadratic Hitchin System}&\longleftarrow&\mbox{Isomonodromy
Equations}\\ \uparrow&&\uparrow\\ \mbox{Quadratic Beilinson--Drinfeld
System}&\longleftarrow&\mbox{Heat Operators/KZ Equations}
\end{array}
\end{equation*}
Namely, the isospectral flows of the Hitchin system arise as a
degeneration of the isomonodromy equations as $\la\to 0$ (see also
\cite{Krich}), and the isomonodromy Hamiltonians appear as a
degeneration of the KZB equations, previously proved for the
Schlesinger equations (isomonodromy for regular singularities in genus
zero) in \cite{Res} (see also \cite{Har}).

\subsubsection{Outline of the paper.}
We begin in \secref{patterns} with a review of the formalism of
localization and the description of twisted differential operators and
symbols from vacuum modules following \cite{Jantzen}. We then apply
this formalism in the special case when the Lie algebras are the
affine Kac-Moody and the Virasoro algebras. First, we review the
necessary features of these Lie algebras in \secref{Vir-KM
algebras}. \thmref{Sugawara embedding} provides the properties of the
Segal-Sugawara construction at the level of representations. Next, we
describe in \secref{moduli spaces} the spaces on which
representations of these Lie algebras localize. These are the moduli
stacks of curves and bundles on curves. In \secref{Sugawara section}
we describe the implications of the Segal-Sugawara constructions for
the sheaves of differential operators on these moduli spaces. Finally,
we consider the classical limits of the localization and their
applications to the isomonodromy questions in \secref{classical limits}.

\medskip

\noindent {\bf Acknowledgments:}\, We would like to thank Tony Pantev,
Spencer Bloch, Sasha Beilinson, Dennis Gaitsgory and Matthew Emerton
for useful discussions, and Roman Fedorov for comments on an early
draft.

\section{Localization.}\label{patterns}

In this section we review the localization construction for
representations of Harish--Chandra pairs as modules over algebras of
twisted differential operators or sheaves on twisted cotangent
bundles. In particular we emphasize the localization of vacuum
representations, which give the sheaves of twisted differential
operators or twisted symbols (functions on twisted cotangent bundles)
themselves.

\subsection{Deformations and limits.}\label{limits}

We first review the standard pattern of deformations and
limits of algebras (see \cite{Jantzen}).

First recall that by the Rees construction, a filtered vector space
$B=\bigcup_i B_i$ has a canonical deformation to its associated graded
$\on{gr}B$ over the affine line ${\mbb A}^1=\on{Spec}\C[\lambda]$,
with the fiber at the point $\la$ isomorphic to $B$ for $\la\neq 0$,
and to $\on{gr} B$ for $\la=0$. The corresponding $\C[\la]$--module is
just the direct sum $\bigoplus_i B_i$, on which $\la$ acts by mapping
each $B_i$ to $B_{i+1}$.

Let $\g$ be any Lie algebra. The universal enveloping algebra $U\g$
has the canonical PBW filtration. The Rees construction then gives us
a one-parameter family of algebras $U_\la\g$, with the associated
graded algebra $U_0\g=\on{Sym}\g$ being the symmetric algebra. The
$\la$--deformation rescales the bracket by $\la$, so that taking the
$\la$--linear terms defines the standard Poisson bracket on the
$\Sym\g=\C[\g^*]$.

Let
$$0\to\C{\bf 1}\to \ghat\to \g\to 0$$ be a one--dimensional central
extension of $\g$. We define a two--parameter family of algebras
$$U_{k,\la}\ghat=U_\la\ghat/({\bf 1}-k\cdot 1),$$ which specializes to
$U_\la \g$ when $k=0$. When $\la=1$ we obtain the family $U_k\ghat$ of
level $k$ enveloping algebras of $\ghat$, whose representations are
the same as representations of $\ghat $ on which $\bf 1$ acts by
$k \cdot \on{Id}$.  This family may be extended to
$(k,\la)\in\Pone\times\Aone$, by defining a $\C[k,k\inv]$--algebra
$$\wh{U}=\oplus \left(\frac{\la}{k}\right)^i (U_{k,\la}\ghat)_{\leq
i}$$ and setting
$$U_{\infty,\la}\ghat=\wh{U}/k\inv\wh{U}.$$ Equivalently, we consider
the $\C[k\inv]$--lattice $\Lambda$ inside $U_{k,\la}$ (for
$k\in\Aone\sm 0$) generated by $1$ and $\ol{x}=\frac{\la}{k}x$ for
$x\in\g$. The algebra $U_{\infty,\la}\ghat$ is then identified with
$\Lambda/k\inv\Lambda$.  (Note that the algebras $U_{k,\la}\ghat$ may
also be obtained from the standard $\la$--deformation of the filtered
algebras $U_k\ghat=U_{k,1}\ghat$.)

Let us choose a vector space splitting $\ghat\cong\g\oplus \C{\bf 1}$,
so that the bracket in $\ghat$ is given by a two--cocycle $\langle
\cdot,\cdot\rangle$ on $\g$:
$$[x,y]_{\ghat}=[x,y]_\g+\langle x,y\rangle\cdot{\bf 1}$$ for
$x,y\in\g\subset\ghat$.  If we let $\ol{x}$ denote $\frac{\la}{k}x$
for $x\in\g\subset\ghat$, the algebra $U_{k,\la}\ghat$ may be
described as generated by elements $\ol{x}, x \in \g$ and $\bf 1$,
with the relations
$$[\ol{x},\ol{y}]=\frac{\la}{k}(\ol{[x,y]_\g}+\la \langle
x,y\rangle\cdot{\bf 1}), \qquad [{\bf 1},\ol{x}] = 0.$$
In particular, in the limit $k \to \infty$ we obtain a commutative
algebra
$$U_{\infty,\la}\ghat\cong \on{Sym}_{\la}\ghat= \on{Sym}\ghat/({\bf
1}-\la)$$ with a Poisson structure.

The algebra $\on{Sym}_{\la}\ghat$ is nothing but the algebra of
functions $\C[\ghat^*_{\la}]$ on the hyperplane in $\ghat^*$
consisting of functionals, whose value on $\bf 1$ is $\la$. The
functions on $\ghat^*$ form a Poisson algebra, with the
Kirillov--Kostant Poisson bracket determined via the Leibniz rule by
the Lie bracket on the linear functionals, which are elements of
$\ghat$ itself. Since $\bf 1$ is a central element,
$\on{Sym}_{\la}\ghat$ inherits a Poisson bracket as well.

Thus, we have defined a two--parameter family of
``twisted enveloping algebras'' $U_{k,\la}\ghat$ for
$(k,\la)\in\Pone\times \Aone$ specializing to $U_k\ghat$ for $\la=1$
and to the Poisson algebras $\on{Sym}_\la\ghat$ for $k=\infty$.
The algebras $U_{k,0}\ghat$ with $\la=0$ are all isomorphic (as
Poisson algebras) to the symmetric algebra $\Sym\g$.

\subsubsection{Harish-Chandra pairs.}\label{HC pairs}
We now suppose that the Lie algebra $\g$ is part of a Harish-Chandra
pair $(\g,K)$.  Thus $K$ is an algebraic group, and we are given an
embedding $\kk=\on{Lie}K\subset \g$ and an action of $K$ on $\g$
compatible with the adjoint action of $K$ on $\kk$ and the action of
$\kk$ on $\g$.  Suppose further that we have a central extension
$\ghat$, split over $\kk$. Then the action of $K$ lifts to $\ghat$ and
$(\ghat,K)$ is also a Harish-Chandra pair.

It follows that the subalgebra of $U_{k,\la}\ghat$ generated by the
elements $\ol{x}$ for $x\in\kk$ is isomorphic to $U_{\la/k}\kk$,
degenerating to $\Sym\kk$ for $\la=0$ or $\la=\infty$. Moreover, the
action of $\kk$ on $U_{k,\la}\ghat$ given by bracket with
$x=\frac{k}{\la}\ol{x}$ preserves the lattice $\Lambda$ and thus is
well-defined in the limit $k=\infty$. (In this limit the action is
generated by Poisson bracket with
$\kk\subset\Sym\kk\subset\Sym_{\la}\ghat$.) In particular $\kk$ acts
on $U_{k,\la}\ghat$ for all $(k,\la)\in\Pone\times\Aone$, and so it
makes sense to speak of $(U_{k,\la}\ghat,K)$--modules.

\subsection{Twisted differential operators and twisted cotangent
bundles.}\label{TDOs}

We review briefly the notions of sheaves of twisted differential
operators and of twisted symbols (functions on twisted cotangent
bundles) following \cite{Jantzen}.

Let $M$ be a smooth variety equipped with a line bundle $\Ll$.
The Atiyah sequence of $\Ll$ is an extension of Lie algebras
$$0\to \Oo_M\to \T_\Ll\stackrel{a}{\longrightarrow} \Theta_M\to 0,$$
where $\T_\Ll$ is the Lie algebroid of infinitesimal automorphisms of $\Ll$,
$\Theta_M$ is the tangent sheaf, and the anchor map $a$ describes the
action of $\T_\Ll$ on its subsheaf $\Oo_M$.

We can generalize the construction of the family of algebras
$U_{k,\la}\ghat$ from \secref{limits} to the Lie algebroid $\T_\Ll$.
We first introduce the sheaf of unital associative algebras $\D_k$
generated by the Lie algebra $\T_\Ll$, with $1\in\Oo_M\subset\T_\Ll$
identified with $k$ times the unit. Due to this identification, $\D_k$
is naturally a filtered algebra. Thus we have the Rees
$\la$--deformation $\D_{k,\la}$ from $\D_k$ to its associated graded
algebra at $\la=0$.  Specializing to $\lambda=0$ we obtain the
symmetric algebra $\Sym\Theta_M$ with its usual Poisson structure,
independently of $k$.  We may also define a limit of the algebras
$\D_{k,\la}$ as $k$ goes to infinity.  In order to do so we introduce
a $\C[k,k\inv]$--algebra $\D_{\la,k}$ as $$\D_{\la,k} =
\bigoplus_{i=0}^\infty \left(\dfrac{\la}{k}\right)^i\cdot \D_{\leq
i}(\Ll^{\otimes k}),$$ and define the limit algebra
$$\D_{\infty,\la}=\D_{\la,k}/ k\inv\cdot \D_{\la,k}.$$ This algebra is
commutative, and hence inherits a Poisson structure from the
deformation process.

\subsubsection{Twisted differential operators.}
Consider the sheaf $\D(\Ll)$ of differential operators acting on
sections of $\Ll$. This is a sheaf of filtered associative algebras,
and the associated graded sheaf is identified with the sheaf of
symbols (functions on the cotangent bundle) as Poisson
algebras. Similarly, the sheaf of differential operators
$\D(\Ll^{\otimes k})$ for any integer $k$ is also a sheaf of filtered
associative unital $\Oo$--algebras, whose associated graded is
commutative and isomorphic (as a Poisson algebra) to the sheaf of
functions on $T^*M$.  The sheaf $\D_k$ defined above for an arbitrary
$k\in \C$ shares these properties, which describe it as a {\em sheaf
of twisted differential operators}. For $k\in\Z$ there is a canonical
isomorphism of sheaves of twisted differential operators $\D_k\cong
\D(\Ll^{\ot k})$.  In particular, the subsheaf $\D(\Ll^{\ot k})_{\leq
1}$ of differential operators of order at most one is identified with
the subsheaf $(\D_k)_{\leq 1}$, which is nothing but $\T_{\Ll}$ with
the extension rescaled by $k$. The sheaf $\D_k$ for general $k$ may
also be described as a subquotient of the sheaf of differential
operators on the total space of $\Ll^{\times}$, the complement of the
zero section in $\Ll$ (namely, as a quantum reduction by the natural
action of $\C^{\times}$).

\subsubsection{Twisted cotangent bundles.}
We may also associate to $\Ll$ a {\em twisted cotangent bundle} of
$M$. This is the affine bundle $\Conn\Ll$ of connections on $\Ll$,
which is a torsor for the cotangent bundle $T^*M$.  For $\la\in\C$ we
have a family of twisted cotangent bundles $T^*_\la M$, which for
$\la\in\Z$ are identified with the $T^*M$--torsors $\Conn\Ll^{\otimes
\la}$ of connections on $\Ll^{\otimes \la}$. Let $\Ll^{\times}$ denote
the principal $\Cx$--bundle associated to the line bundle $\Ll$. Then
the space $T^*_\la M$ is identified with the Hamiltonian reduction (at
$\la\in\C=(\on{Lie}\Cx)^*$) of the cotangent bundle $T^*\Ll^{\times}$
by the action of $\Cx$. In particular, $T^*_\la M$ carries a canonical
symplectic structure for any $\la$.

Thus the {\em twisted symbols}, i.e., functions on $T^*_\la M$ (pushed
forward to $M$) form a sheaf of $\Oo_M$--Poisson algebras $\Oo(T^*_\la
M)$. Its subsheaf of {\em affine} functions $\Oo(T^*_\la M)_{\leq 1}$
on the affine bundle $T^*_\la M$ forms a Lie algebra under the Poisson
bracket.  It is easy to check that there is an isomorphism of sheaves
of Poisson algebras $\D_{\infty,\la}\cong \Oo(T^*_\la M)$, identifying
$\Oo(T^*_\la M)_{\leq 1}$ with $\T_{\Ll}$ (with the bracket rescaled
by $\la$).

\bigskip

Thus to a line bundle we have associated a two--parameter family of
algebras $\D_{k,\la}$ for $(k,\la)\in \Pone\times \Aone$.  This
specializes for $k\in\Z,\la=1$ to the differential operators
$\D(\Ll^{\ot k})$, for $\la=0$ and $k$ arbitrary to symbols
$\Oo(T^*M)$, and for $k=\infty, \la\in\Z$ to the Poisson algebra of
twisted symbols, $\Oo(\Conn\Ll^{\ot \la})$.

\subsection{Localization.}\label{localization}

In this section we combine the algebraic picture of \secref{limits}
with the geometric picture of \secref{TDOs} following \cite{Jantzen,
Hecke} (see also \cite{book}, Ch. 16).

Let $\g,\ghat,K$ be as in \secref{limits}.  We will consider the
following geometric situation: $\Mhat$ is a smooth scheme equipped
with an action of $(\g,K)$, in other words, an action of the Lie
algebra $\g$ on $\Mhat$ which integrates to an algebraic action of the
algebraic group $K$. Let $\M=\Mhat/K$ be the quotient, which is a
smooth algebraic stack, so that $\pi:\Mhat\to \M$ is a $K$--torsor.
Thus if the $\g$ action on $\Mhat$ is infinitesimally transitive, then
given $\wh{x}\in \Mhat$ and $x\in \M$ its $K$--orbit, the formal
completion of $\M$ at $x$ is identified with the double quotient of
the formal group $\on{exp}(\g)$ of $\g$: by the formal group
$\exp(\kk)$ of $\kk$ on one side, and the formal stabilizer of
$\wh{x}$ on the other.

Let $\Ll$ denote a $K$--equivariant line bundle on $\Mhat$ and
the corresponding line bundle on $\M$.  We assume that the action of
$\g$ on $\Mhat$ lifts to an action of $\ghat$ on $\Ll$, so that ${\bf
1}$ acts as the identity.  For $k,\la\in\Pone\times\Aone$, we have an
algebra $U_{k,\la}\ghat$ from \secref{limits} and a sheaf of algebras
$\D_{k,\la}$ on $\Mhat$ from \secref{TDOs}. The action of $\ghat$ on
$\Ll$ matches up these definitions: we have an algebra homomorphism
$$\Oo_{\Mhat}\ot U_{k,\la}\ghat\longrightarrow \D_{k,\la},\hskip.3in
(k,\la)\in\Pone\times\Aone.$$ This may also be explained by
reduction. Namely, $\ghat$ acts on the total space $\Ll^{\times}$ of
the $\Cx$--bundle associated to $\Ll$.  Therefore $U\ghat$ maps to
differential operators on $\Ll^{\times}$. It also follows that $\ghat$
acts in the Hamiltonian fashion on $T^*\Ll^{\times}$, so that
$\Sym\ghat$ maps in the Poisson fashion to symbols on
$\Ll^{\times}$. The actions of $U_k\ghat\to \D_k$ and
$\Sym_\la\ghat\to \Oo(T^*_\la\Mhat)$ are then obtained from the
quantum (respectively, classical) Hamiltonian reduction with respect
to the $\Cx$ action on $\Ll^{\times}$ and the moment values $k$
(respectively, $\la$).

\subsubsection{Localization functor}

The action of the Harish-Chandra pair $(\ghat,K)$ on $\Mhat$ and $\Ll$
may be used to construct localization functors from
$(\ghat,K)$--modules to sheaves on $\M$. Let $M$ be a
$U_{k,\la}\ghat$--module with a compatible action of $K$. We then
consider the sheaf of $\D_{k,\la}$--modules
$$\wh{\Delta}(M) = \D_{k,\la} \und{U_{k,\la}\ghat}\ot M$$ on $\Mhat$,
whose fibers are the coinvariants of $M$ under the stabilizers of the
$\ghat$--action on $\wh{\Ll}^\times$. (By our hypotheses, these lift
to $\ghat$ the stabilizers of the $\g$--action on $\Mhat$.) The sheaf
$\wh{\Delta}(M)$ is a $K$--equivariant $\D_{k,\la}$--module, and so it
descends to a $\D_{k,\la}$--module $\Delta_{k,\la}(M)$ on $\M$, the
localization of $M$:
$$\Delta_{k,\la}(M)=(\pi_*(\D_{k,\la} \und{U_{k,\la}\ghat}\otimes
M))^K.$$

Equivalently, we may work directly on $\M$ by twisting.  Let
 $$\Mcal=\Mhat \und{K}\times M=\pi_*(\Oo_{\Mhat}\ot M)^K$$ denote the
twist of $M$ by the $K$--torsor $\Mhat$ over $\M$ (i.e., the vector
bundle on $\M$ associated to the principal bundle $\Mhat$ and
representation $M$). Then there is a surjection
\begin{equation}\label{twist to coinvariants}
\Mcal\twoheadrightarrow \Delta_{k,\la}(M),
\end{equation}
and $\Delta_{k,\la}(M)$ is the quotient of $\Mcal$ by the twist
$(\g_{\on{stab}})_{\Mhat}\subset (\g)_{\Mhat}$ of the sheaf of
stabilizers.


\subsubsection{Remark.}
The localization functor $\Delta_{k,\la}$ ``interpolates'' between the
localization $$\Delta_k: (U_k\ghat,K)-\on{mod} \quad \longrightarrow
\quad \D_k-\on{mod}$$ as modules for twisted differential operators, and
the classical localization
$$\ol{\Delta}_\la: (\Sym_\la\ghat,K)-\on{mod} \quad \longrightarrow
\quad \Oo(T^*_\la \M)-\on{mod}$$ as quasicoherent sheaves on $T^*_\la
\M$.  The latter assigns to a module over the commutative ring
$\Sym_\la\ghat$ a quasicoherent sheaf on $T^*_\la\Mhat$ via the
embedding of $\Sym_{\la}\ghat$ into global twisted symbols. This
$K$--equivariant sheaf descends to $\M$, where it becomes a module
over the sheaf of twisted symbols $\Oo(T^*_\la \M)$ and therefore we
obtain a sheaf on the twisted cotangent bundle.

\subsection{Localizing the vacuum module.}\label{localizing vacuum}

The fundamental example of a $(\g,K)$--module is the {\em vacuum
module} $$\Vgk=\on{Ind}_{U\kk}^{U\g}\C=U\g/(U\g\cdot \kk)$$ (we will
denote it simply by $V$ when the relevant Harish-Chandra pair $(\g,K)$
is clear). The vector in $\Vgk$ corresponding to $ 1\in\C$ is denoted
by $\vac$ and referred to as the vacuum vector. It is a cyclic vector
for the action of $\g$ on $\Vgk$.  The vacuum representation has the
following universality property: given a $\kk$--invariant vector $m\in
M$ in a representation of $\g$, there exists a unique
$\g$--homomorphism $m_*:\Vgk\to M$ with $m_*(\vac)=m$.

An important feature of the localization construction for
Harish-Chandra modules is that the sheaf of twisted differential
operators $\D_k$ itself arises as the localization $\Delta(\Vgk)$ of
the vacuum module.  This is a generalization of the description of
differential operators on a homogeneous space $G/K$ as sections of the
twist
$${\mc V}_{\g,K}=\Pc_K \und{K}\times \Vgk$$
by the $K$--torsor $\Pc_K=G$ over $G/K$.
Informally, the general statement follows by applying reduction from
$G/K$ to $\M$, which is (formally) a quotient of the formal homogeneous
space $\on{exp}(\g)/\on{exp}(\kk)$.

Recall that $\pi:\Mhat\to \M$ denotes a $K$--torsor with compatible
$\g$--action. We assume $K$ is connected and affine. In particular,
the functor $\pi_*$ is exact.

\subsubsection{\bf Proposition.} \label{twists are algebras}
{\em
\begin{enumerate}
\item The twist $\V=\Mhat \und{K}\times V$, its subsheaf of invariants
  $\V^K=\Oo_\M\ot V^K$, and
  the localization $\Delta(V)$ are sheaves of algebras on $\M$.
\item The localization $\Delta(V)$ is naturally
  identified with the sheaf $\D$ of differential operators on $\M$.
\item The natural maps
$$\Oo_\M\ot V^K \longrightarrow \V\longrightarrow \Delta(V)=\D_\M$$
are algebra homomorphisms.
\item For any $(\g,K)$--module $M$, the sheaves
$$\Oo_\M\ot M^K\longrightarrow {\Mcal}\longrightarrow \Delta(M)$$
are modules over the corresponding algebras in $(3)$.
\end{enumerate}
}

\subsubsection{Proof.}
The vector space $V^K=[U\g/U\g\cdot\kk]^\kk$ is naturally an algebra,
isomorphic to the quotient $N(I)/I$ of the normalizer $N(I)$ of the
ideal $I=U\g\cdot \kk$ by $I$. It is also isomorphic to the algebra
$(\on{End}_{\g}(V))^{\on{opp}}$ (with the opposite multiplication),
thanks to the universality property of $V$, and thus to the algebra of
endomorphisms of the functor of $K$--invariants (with the opposite
multiplication). Hence $\Oo_\M\ot V^K$ is an $\Oo_\M$--algebra.

The sheaf $\Oo_{\Mhat}\ot U\g$ is a sheaf of algebras, where $U\g$ acts
on $\Oo_{\Mhat}$ as differential operators, and so is ${\mc
  A}=\pi_*(\Oo_{\Mhat}\ot U\g)$. Let ${\mc K}$ be
the sheaf of $\Oo_\M$--Lie algebras
$\pi_*(\Oo_{\Mhat}\ot\kk)$ (note that the action of $\kk$ on $\Oo_{\Mhat}$ is
$\pi\inv(\Oo_\M)$--linear). Then we have
\begin{eqnarray*}
\V&=&[\pi_*(\Oo_{\Mhat}\ot V)]^K\\ &=&[\pi_*(\Oo_{\Mhat}\ot V)]^{\mc
K}\\ &=&[{\mc A}/{\mc A}\cdot {\mc K}]^{\mc K},
\end{eqnarray*}
which is thus a sheaf of algebras (again as the normalizer of an ideal modulo that ideal).

On the other hand, the localization
\begin{eqnarray*}
\Delta(V) &=& [\pi_*(\D_{\Mhat}\ot_{U\g} V)]^K\\
&=& [\pi_*(\D_{\Mhat}/ \D_{\Mhat} \cdot{\mc K})]^{\mc K}\\
&=& \D_\M
\end{eqnarray*}
since $\mc K$ is the sheaf of {\em vertical} vector fields. This is in
fact the common description of differential operators on a quotient as
the quantum hamiltonian (BRST) reduction of differential operators
upstairs.  The obvious maps are compatible with this definition of
the algebra structures. This proves parts (1)--(3) of the proposition.

Finally, for part (4), observe that the $U\g$--action on $M$ descends
to a $\V$--action on $\mc M$, since the ideal $\mc K$ acts trivially
on the $\mc K$--invariants ${\mc M}=\pi_*(\Oo_{\Mhat}\ot M)^{\mc K}$.
Equivalently, the Lie algebroid $\Mhat \und{K}\times \g$ on $\M$ acts
on the twist $\mc M$, and the action descends to the quotient
algebroid $\Mhat \und{K}\times \g/\kk$, which generates $\V$. The
compatibility with the actions of invariants on ${\mc M}^K$ and of the
localization $\D$ on the $\D$--module $\Delta(W)$ follow.

\subsubsection{\bf Corollary.}    \label{global diff}
{\em The space of invariants $V^K$ is naturally a subalgebra of the
algebra $\Gamma(\M,\D)$ of global differential operators on $\M$.}

\subsubsection{Changing the vacuum.}

Let $\ol{K}\subset K$
be a normal subgroup.  Then a minor variation of the proof establishes
that the twist of $\ol{K}$--invariants $\V^{\ol{K}}\hookrightarrow \V$
is a subalgebra.  Note moreover that $\V^{\ol{K}}$ depends only on the
induced $K/\ol{K}$--torsor ${\wh{\mf N}}=(K/\ol{K})_{\Mhat}$:
$$\V^{\ol{K}}=\wh{\mf N}\times_{K/\ol{K}} V^{\ol{K}}.$$

Now suppose $(\ab,K)$ is a sub--Harish-Chandra pair of $(\g,K)$ (i.e.
$\kk\subset\ab\subset \g$). Then we have the two vacuum modules
$V=V_{\g,K}\supset V_{\ab,K}$ and the module
$V_{\g,\ab}=\on{Ind}_{\ab}^{\g}\C$, which is not literally the vacuum
module of a Harish-Chandra pair unless $\ab$ integrates to an
algebraic group.

By \propref{twists are algebras}, the twist $\V_{\ab,K}$ and the
localization $\Delta(V_{\ab,K})$ -- with respect to the
Harish-Chandra pair $(\ab,K)$ -- are sheaves of algebras, and act on
$\V$ and $\Delta(V)$ respectively. We then have

\subsubsection{\bf Lemma.}\label{change of vacuum}
{\em The quotients $\V/(\V_{\ab,K}\cdot \V)$ and
$\Delta(V)/(\Delta(V_{\ab,K}) \cdot \Delta(V))$ are isomorphic,
respectively, to the twist $\V_{\g,\ab}$ and the localization
$\Delta(V_{\g,\ab})$ with respect to the Harish-Chandra pair
$(\g,K)$.}

\subsubsection{Twisted and deformed vacua.}

In general, we define the twisted, classical and deformed vacuum
representations $V_k(\ghat)$, $\ol{V}_{\la}(\ghat)$ and $V_{k,\la}$ of
$U_k\ghat$, $\Sym_\la\ghat$ and $U_{k,\la}\ghat$, respectively, as the
induced representations
\begin{eqnarray*}
V_{k}&=&U_k\ghat \und{U\kk}\ot \C \\
\ol{V}_\la&=& \Sym_\la\ghat \und{\Sym\kk}\ot \C \\
V_{k,\la}&=& U_{k,\la}\ghat \und{U_{k,\la}\khat}\ot \C,
\end{eqnarray*}
where we use the observation (\secref{HC pairs}) that $U_{k}\khat$ is
isomorphic to $U\kk$ for $k\neq\infty$ and to $\Sym\kk$ for
$k=\infty$. All of these representations carry compatible
$K$--actions.

We then have a general localization principle, describing twisted
differential operators and twisted symbols on $\M$ in terms of the
corresponding vacuum representations. The proof of \propref{twists are
algebras} generalizes immediately to the family $V_{k,\la}$ over
$\Pone\times\Aone$. This is deduced formally from the analogous
statement for homogeneous spaces $G/K$, twisted by an equivariant line
bundle. In that case the algebra of functions on the cotangent space
$T^*G/K|_{1K}=(\g/\kk)^*$ is identified with
$$\ol{V}_0=\on{Sym}(\g/\kk)= \on{Sym}\g \und{\on{Sym}\kk}\ot \C,$$ and
sections of the twist of $\ol{V}_0$ give the sheaf of functions on
$T^*G/K$.  Note that $\ol{V}_\la$ is a (commutative) algebra; however,
it only becomes Poisson after twisting or localization. To summarize,
we obtain the following

\subsubsection{\bf Proposition.}\label{localization of vacuum in general}

{\em The localization $\Delta_{k,\la}(V_{k,\la})$ is canonically
isomorphic to the sheaf of algebras $\D_{k,\la}$, for
$(k,\la)\in\Pone\times\Aone$. In particular
$\Delta_k(V_{k})\cong\D_k$, and $\ol{\Delta}_{\la}(\ol{V}_{\la})\cong
\Oo(T^*_\la \M)$ as Poisson algebras.}

\section{Virasoro--Kac-Moody Algebras.}\label{Vir-KM algebras}

In this section we introduce the Lie algebras to which we wish to
apply the formalism of localization outlined in the previous
sections. These are the affine Kac-Moody algebras, the Virasoro
algebra and their semi-direct product. We describe the Segal-Sugawara
construction which expresses the action of the Virasoro algebra on an
affine algebra as an ``internal'' action. We interpret this
construction in terms of a homomorphism between vacuum representations
of the Virasoro and Kac-Moody algebras, and identify the critical and
classical limits of these homomorphisms.  In the subsequent sections
we will use the localization of this construction to describe sheaves
of differential operators on the moduli spaces of curves and bundles
on curves.

\subsection{Virasoro and Kac-Moody algebras.}

Let $\g$ be a simple Lie algebra. It carries an invariant bilinear
form $(\cdot,\cdot)$ normalized in the standard way so that the square
length of the maximal root is equal to 2. We choose bases
$\{J^a\},\{J_a\}$ dual with respect to this bilinear form. The affine
Kac-Moody Lie algebra $\ghat$ is a central extension
$$0\to \C{\bf K} \to \ghat\to L\g\to 0$$ of the loop algebra
$L\g=\g\ot\K$ of $\g$, with (topological) generators $\{J^a_n=J^a\ot
t^n\}_{n\in\Z}$ and relations
\begin{equation}\label{KM relations}
[J^a_n,J^b_m]=[J^a,J^b]_{n+m}+(J^a,J^b)\delta_{n,-m}{\bf K}.
\end{equation}
The central extension splits over the Lie subalgebra $\g(\Oo)\subset
\g(\K)$ (topologically spanned by $J^a_n, n\geq 0$), so that the
affine proalgebraic group $G(\Oo)$ acts on $\ghat$.  Thus we have a
Harish-Chandra pair $(\ghat, G(\Oo))$.  We denote by
$$V_k=V_k(\ghat)=U_k\ghat \und{U\g(\Oo)}\ot \C$$ the corresponding
vacuum module (see \secref{localizing vacuum}).  More generally, let
$\mm\subset \Oo$ denote the maximal ideal $t\C[[t]]$. Let
$\g_n(\Oo)=\g\ot\mm^n\subset\g(\Oo)$ denote the congruence subalgebra
of level $n$, and $\Gn\subset G(\Oo)$ the corresponding algebraic
group, consisting of loops which equal the identity to order $n$.  We
denote by $$V_k^n=V_k^n(\ghat)=U_k\ghat \und{U\g_n(\Oo)}\ot \C$$ the
vacuum module for $(\ghat,\Gn)$.

Let $Vir$ denote the Virasoro Lie algebra. This is a central extension
$$0\to \C{\bf C}\to Vir\to \on{Der}\K \to 0$$ of the Lie algebra of
derivations of the field $\K$ of Laurent series. It has (topological)
generators $L_n=-t^{n+1}\pa_t, n\in\Z$, and the central element ${\bf
C}$, with relations
\begin{equation}\label{Vir relations}
[L_n,L_m]=(n-m)L_{n+m}+\frac{1}{12}(n^3-n)\delta_{n,-m}{\bf C}.
\end{equation}
The central extension splits over the Lie subalgebra
$\on{Der}\Oo\subset\on{Der}\K$, topologically spanned by $L_n, n\geq
-1$. Consider the module
$$Vir_c=U_cVir \und{U\on{Der}\Oo}\ot \C.$$ Since $\DerO$ is not the Lie
algebra of an affine group scheme, $Vir_c$ is not strictly speaking a
vacuum module of a Harish--Chandra pair.

The affine group scheme $\AutO$ of all changes of coordinates on the
disc fixes the closed point $t=0$, so that $L_{-1}=-\pa_t$ is not in
its Lie algebra, which we denote by $\Der_1\Oo$.  More generally, we
set
$$\Der_n(\Oo)=\mm^{n}\Der\Oo\cong t^{n}\C[[t]]\partial_t$$ and let
$\Aut_n(\Oo)\subset \Aut_0(\Oo)=\AutO$ be the corresponding algebraic
subgroups (for $n\geq 0$).  Then we have Harish-Chandra pairs
$(Vir,\Aut_n\Oo)$ and the corresponding vacuum modules $Vir_c^n$ (thus
we may denote $Vir_c=Vir_c^{0}$).

The Virasoro algebra acts on $\ghat$, via its action as derivations of
$\K$. In terms of our chosen bases this action is written as follows:
\begin{equation*}
[L_n,J^a_m]=-mJ^a_{n+m}.
\end{equation*}
Let $\gtil$ be the resulting semidirect product Lie algebra $Vir
\ltimes \ghat$. Thus $\gtil$ has topological generators
$$J^a_m\hskip.1in(m\in\Z,J^a\in\g),\hskip.3in L_n\hskip.1in
(n\in\Z),\hskip.3in {\bf K},{\bf C},$$
with relations as above.

The central extension splits over the semidirect product
$\gtil_+=\on{Der}\Oo\ltimes \g(\Oo)$.  We let
$$V_{k,c}(\gtil)=\on{Ind}_{\gtil_+\oplus\C{\bf 1}}^{\gtil}\C_{k,c}
=U\gtil \und{U(\gtil_+\oplus\C{\bf 1})}\otimes \C_{k,c},$$ where
$\C_{k,c}$ is the one--dimensional representation of
$\gtil_+\oplus\C{\bf K}\oplus\C{\bf C}\subset\gtil$ on which
$\gtil_+$ acts by zero and ${\bf K,C }$ act by $k,c$.  Let
$$\gtil_n(\Oo)=\Der_{2n}(\Oo)\ltimes\g_n(\Oo)=\mm^{2n}\DerO\ltimes
\mm^n\g(\Oo).$$

\subsubsection{Remark.} This normalization above, by which
we pair $\g_n(\Oo)$ with $\Der_{2n}(\Oo)$, is motivated by the
Segal-Sugawara construction, cf. \propref{quadratic level}: the
Virasoro generators will be constructed from {\em quadratic}
expressions in $\ghat$. Thus in geometric applications of the
Segal--Sugawara operators the orders of trivializations or poles along
the Virasoro (moduli of curves) directions will turn out to be double
those along the Kac--Moody (moduli of bundles) directions (for example
the quadratic Hitchin hamiltonians double the order of pole of a Higgs
field).

\subsubsection{}The corresponding representation
$$V_{k,c}^n=U_{k,c}\gtil \und{U\gtil_n(\Oo)}\ot\C$$ agrees with
$V_{k,c}$ for $n=0$, while for $n>0$ it is identified as the vacuum
module for the Harish-Chandra pair $(\gtil,\Gtil_n(\Oo))$ where
$\Gtil_n(\Oo)=\Aut_{2n}\Oo\ltimes\Gn$. We will also denote by
$\Gtil(\Oo)$ the semidirect product $\AutO\ltimes G(\Oo)$.

The pattern of deformations and limits from \secref{limits} applies to
the Kac-Moody and Virasoro central extensions, giving rise to two
families of algebras which we denote $U_{k,\la}\ghat$ and
$U_{c,\mu}Vir$. Thus for $\la=\mu=1$ and $k=c=0$,
$U_{k,\la}\ghat=U\g(\K)$ and $U_{c,\mu}Vir=U\on{Der}\K$, etc. For
$k=c=\infty$ we obtain the Poisson algebras $\on{Sym}_{\la}\ghat$ and
$\on{Sym}_{\mu}Vir$.

We also have classical vacuum modules $\ol{V}_\la(\ghat)$ for
$(\Sym\ghat,G(\Oo))$ and $\ol{V}_\mu(Vir)$ for $(\Sym Vir,\AutO)$ and
the interpolating families $V_{k,\la}(\ghat)$ and $V_{c,\mu}(Vir)$.


\subsection{Vertex algebras.}
The vacuum representations $V_k,Vir_c$ and $V_{k,c}$ have natural
structures of vertex algebras, and $V_k^n,Vir_c^n$ and $V_{k,c}^n$ are
modules over the respective vertex algebras (see \cite{book} for the
definition of vertex algebras and in particular the Virasoro and
Kac-Moody vertex algebras).

The vacuum vector $\vac$ plays the role of the unit for these vertex
algebras. Moreover, all three vacuum modules carry natural
Harish-Chandra actions of $(\Der\Oo,\AutO)$: these are defined by the
natural action of $\Der\Oo$ on $\ghat,Vir$ and $\gtil$, preserving the
positive halves and hence giving rise to actions on the corresponding
representations. These actions give rise to a grading operator $L_0$
and a translation operator $T=L_{-1}=-\partial_z \in\on{Der}\Oo$. The
vertex algebra $V_k$ is generated by the fields
$$J^a(z)=Y(J^a_{-1}\vac,z)=\sum_{n\in\Z} J^a_n z^{-n-1}$$
associated to the vectors $J^a_{-1}\vac\in V_k$, which satisfy the
operator product expansions
$$J^a(z) J^b(w) = \frac{k(J^a,J^b)}{(z-w)^2} +
\frac{[J^a,J^b](w)}{z-w} + \cdots$$ (where the ellipses denote regular
terms).  This may be seen as a shorthand form for the defining
relations \eqref{KM relations}.  This structure extends to the family
$V_{k,\la}$.
Explicitly, introduce
the $\C[\la,k\inv]$--lattice $\Lambda$ in $V_{k,\la}$ generated
by monomials in $\ol{J}_n^a=\frac{\la}{k}J^a_n$.  These satisfy
relations
$$[\ol{J}^a_n,\ol{J}^b_m]=\dfrac{\la}{k} (\ol{[J^a,J^b]}_{n+m}+\la
(J^a,J^b)) \delta_{n,-m},$$
or in shorthand
$$[\ol{J},\ol{J}]=\frac{\la}{k}(\ol{J}+\la(\cdots))$$ The vertex
algebra structure is defined by replacing the $J$ operators by the
rescaled versions $\ol{J}$. In particular we see that we recover the
commutative vertex algebra $\ol{V}_\la(\g)$ when $k=\infty$. Recall
from \cite{book}, Ch. 15, that a {\em commutative} vertex algebra is
essentially the same as a unital commutative algebra with derivation
and grading. Thus the classical vacuum representations
$\ol{V}_\la(\g)$, $\ol{Vir}_{\mu}$ and $\ol{V}_{\la,\mu}(\g)$ are
naturally commutative vertex algebras.

Moreover, the description as a degeneration endows $\ol{V}_\la(\g)$
with a vertex Poisson algebra structure. The Poisson structure comes
with the following relations
$$\{\ol{J}^a_n,\ol{J}^b_m\}=\ol{[J^a,J^b]}_{n+m} +
\la(J^a,J^b)\delta_{n,-m}.$$

Likewise, the vertex algebra $Vir_c$ is generated by the field
\begin{equation}\label{T(z)}
T(z)=Y(L_{-2}\vac,z)= \sum_{n\in\Z} L_n z^{-n-2},
\end{equation}
satisfying the operator product expansion
\begin{equation}\label{Vir OPE}
T(z)T(w)=\frac{c/2}{(z-w)^4}+\frac{2T(w)}{(z-w)^2}+
\frac{\partial_wT(w)}{z-w}+\cdots
\end{equation}
encapsulating the relations \eqref{Vir relations}.
The $\mu$--deformation $Vir_{c,\mu}$ is straightforward.
The commutation relations for the generators
$\ol{L}_n=\frac{\mu}{c}L_n$ read
$$[\ol{L}_n,\ol{L}_m]=\dfrac{\mu}{c}(n-m)\ol{L}_{n+m} +
\dfrac{\mu^2}{c} \frac{1}{12}(n^3-n) \delta_{n,-m},$$
or simply
$$[\ol{L},\ol{L}]=\frac{\mu}{c}(\ol{L}+\mu(\cdots)).$$ In the limit
$c=\infty$ we obtain a vertex Poisson algebra structure on $\ol
V_{\mu}(\g)$.  (As a commutative algebra with derivation it is freely
generated by the single vector $\ol{L}_{-2}\vac$.)  The Poisson
operators $\ol{L}$ have relations
$$\{\ol{L},\ol{L}\}=\ol{L}+\mu(\cdots).$$

The space $V_{k,c}$
also carries a natural vertex algebra structure, so that
$Vir_c\subset V_{k,c}(\g)$ is a vertex subalgebra, complementary to
$V_k(\ghat)\subset V_{k,c}(\g)$, which is a vertex algebra ideal.
In particular, the vertex algebra is generated by the fields $J^a(z)$
and $T(z)$, with the additional relation
\begin{equation}\label{Sugawara OPE}
T(z)J^a(w)=\frac{J^a(w)}{(z-w)^2}+
\frac{\partial_w J^a(w)}{z-w}+\cdots
\end{equation}

We may combine the above deformations into a two--parameter
deformation of $V_{k,c}$. We introduce the $\C[k\inv,c\inv]$--lattice
$\Lambda$ in $V_{k,c}$ generated by monomials in
$\ol{J_n^a}=\dfrac{\la}{k}J^a_n$ and $\ol{L_n}=\dfrac{\mu}{c}L_n$.
These satisfy relations
$$[\ol{J},\ol{J}]=\frac{\la}{k}(\ol{J}+\la),\hskip.2in
[\ol{L},\ol{L}]=\frac{\mu}{c}(\ol{L}+\mu),\hskip.2in
[\ol{L},\ol{J}]=\frac{\mu}{c}\ol{J}.$$


If we impose $\la=\mu=0$, we obtain
the deformation of the enveloping vertex algebra $V_{k,c}$ to the
symmetric vertex Poisson algebra associated to $\wt{\g}$.
We will denote this limit vertex Poisson algebra by $\ol{V}_{0,0}(\g)$.
We may also specialize $\mu$ to $0$, making the $\ol{L}$ generators
central. Thus in this limit we have a noncommutative vertex structure
on the Kac-Moody part, while the Virasoro part degenerates to a
vertex Poisson algebra.  Note however that if we keep $\mu$ nonzero
but specialize $\la=0$, the Kac-Moody generators become commutative
but {\em not} central. Hence the Kac-Moody part does not acquire a
vertex Poisson structure, and we do not obtain an action of the
Kac-Moody vertex Poisson algebra on $V_{k,c}$ in this limit.

In order to obtain a vertex Poisson structure, we need a vertex
subalgebra which becomes central, together with a ``small parameter''
or direction of deformation, in which to take the linear term.  The
``generic limit'' we consider is obtained by letting $k,c\to\infty$,
but with the constraint that our small parameter is
$$\frac{\la}{k}=\frac{\mu}{c}.$$
In other words, we take terms linear in
either one of these ratios. In this limit, the entire vertex algebra
becomes commutative, as can be seen from the relations above. Thus in
the limit we obtain a vertex Poisson algebra, generated by
$\ol{L},\ol{J}$ with relations
$$\{\ol{J},\ol{J}\}=\ol{J}+\la\hskip.2in
\{\ol{L},\ol{L}\}=\ol{L}+\mu\hskip.2in
\{\ol{L},\ol{J}\}=\ol{J}.$$
The resulting vertex Poisson algebra
will be denoted by $\ol V_{\la,\mu}(\g)$, for $\la,\mu\in\C$.

\subsection{The Segal-Sugawara vector.}\label{sug vector}

To any vertex algebra $V$ we associate a Lie algebra $U(V)$
topologically spanned by the Fourier coefficients of vertex operators
from $V$ (see \cite{book}). This Lie algebra acts on any
$V$--module. In the case of the Kac-Moody vertex algebra $V_k$, the
Lie algebra $U(V_k)$ belongs to a completion of the enveloping algebra
$U_k\ghat$. An important fact is that it contains inside it a copy of
the Virasoro algebra (if $k\neq -h^\vee$).  In vertex algebra
terminology, this means that the Kac-Moody vacuum modules are {\em
conformal} vertex algebras.

The conformal structure for $V_{k,c}$ is automatic.  Let
$\omega_V=L_{-2}\vac\in V_{k,c}(\g)$.  This is a conformal vector for
the vertex algebra $V_{k,c}(\g)$: the field $T(z)$ it generates
satisfies the operator product expansion \eqref{Vir OPE}, and hence
its Fourier coefficients $L_n$ (see \eqref{T(z)}) give rise to a
Virasoro action. The operator $L_0$ is the grading operator and
$L_{-1}$ is the translation operator $T$ on $V_{k,c}$. The action of
$\on{Der}\Oo$ induced by the $L_n$ ($n\geq -1$) preserves the
Kac-Moody part $V_k\subset V_{k,c}$, but the negative $L_n$'s take us
out of this subspace.

For level $k$ not equal to minus the dual Coxeter number $\hv$ of
$\g$, the vertex algebra $V_k$ itself carries a Virasoro action, and
in fact a conformal structure, given by the Segal-Sugawara
vector. This means that we have a conformal vector $\omega_S \in
V_k(\ghat)\subset V_{k,c}(\g)$ such that the corresponding field
satisfies the Virasoro operator product expansion \eqref{Vir OPE}.
This conformal vector is given by
\begin{equation}\label{Sugawara vector}
\omega_S=\frac{S_{-2}}{k+h^\vee},\hskip.3in S_{-2}=\frac{1}{2}\sum_a
J^a_{1}J_{a,-1}.
\end{equation}
The corresponding Virasoro algebra has central charge
$\dfrac{k \on{dim}\g}{k+h^\vee}$. In other words,
\begin{equation*}
Y(\omega_S,z)= \sum_{n\in\Z} L_n^S z^{-n-2} \hskip.2in
(L_n^S\in\on{End}V_{k,c}(\g))
\end{equation*}
and we have
\begin{equation}    \label{LnJ}
[L_n^S,J^a_m] = -m J^a_{n+m},
\end{equation}
\begin{equation}    \label{LnS}
[L_n^S,L_m^S]=(n-m)L_{n+m}^S + \frac{1}{12}
\frac{k\on{dim}\g}{k+\hv} (n^3-n)\delta_{n,m}.
\end{equation}

We will continue to take the vector $\omega_V$ as the conformal
vector for $V_{k,c}$, and use the notation $L_n$ for the coefficients
of the corresponding field $T(z)$. Due to the above commutation
relations, $\omega_S$ is {\em not} a conformal vector for $V_{k,c}$,
because $L_0^S$ and $L_{-1}^S$ do not act as the grading and
translation operators on the Virasoro generators.

Note that the commutator $[L_n,J^a_m]$ coincides with the right hand
side of formula \eqref{LnJ}, and therefore $[L_n,L_m^S]$ is given by
the right hand side of formula \eqref{LnS}.

The above construction of the Virasoro algebra action on $V_k$ is
nontrivial from the Lie algebra point of view. Indeed, the operators
$L_n^S$ are given by infinite sums of quadratic expressions in the
$J^a_n$, which are nonetheless well-defined as operators on $V_k$:
$$Y(\omega_S,z)=\frac{1}{2(k+\hv)}\sum_a :J^a(z)J_a(z): \; ,$$
so that
\begin{eqnarray*}
L_n^S&=&\frac{1}{2(k+\hv)} \sum_a\sum_m :J^a_{m}J_{a,n-m}:\\
&=&\frac{1}{2(k+\hv)} \sum_a \left( \sum_{m<0} J^a_m J_{a,n-m}
+\sum_{m\geq 0} J_{a,n-m} J^a_m \right).
\end{eqnarray*}

In fact, since this action is given by a vertex operator in $V_k$, it
follows immediately that any module over the vertex algebra $V_k, k
\neq -h^\vee$, carries a compatible action of the Virasoro algebra.
This may also be expressed using completions of the enveloping
algebras of $\ghat$ and $\gtil$.  These are the completions which act
on all smooth representations, i.e., those in which every vector is
stabilized by a deep enough congruence subgroup $\Gtil_n(\Oo), n\geq
0$.  Namely, we define a completion of $U_{k,c} \gtil$ as the inverse
limit
$$\wh{U}_{k,c}\gtil=\lim_{\longleftarrow} \;
U_{k,c}\gtil/(U_{k,c}\gtil\cdot \gtil_n(\Oo)).$$ This is a complete
topological algebra, since the multiplication on $U_{k,c}(\gtil)$ is
continuous in the topology defined by declaring the left ideals
generated by $\gtil_n(\Oo)$ to be base of open neighborhoods of 0. We
define a completion of $U_k \ghat$ in the same way. These completions
contain the Lie algebras $U(V_{k,c})$ and $U(V_k)$, and in particular
for $k \neq -h^\vee$ they contain the Virasoro algebra generated by
the Segal-Sugawara operators $L_n^S, n \in \Z$. Hence any smooth
representation of $\gtil$ or $\ghat$ of level $k \neq -h^\vee$
inherits a Virasoro action. In particular, the algebra
$\wh{U}_{k,c}(\gtil)$ acts on the vacuum modules $V_{k,c}^n$

\subsubsection{\bf Proposition.}\label{quadratic level}
{\em For any $k \neq -h^\vee$, the Segal-Sugawara operators $L_m^S, m
\in \Z$, define $(Vir,\Aut_{2n}\Oo)$--action of central charge
$c=\frac{k\dim\g}{k+\hv}$ on the vacuum modules $V_k^n$ and
$V_{k,c}^n$. Together with the action of $(\ghat,G_n(\Oo))$, this
$(Vir,\Aut_{2n}\Oo)$--action combines into a
$(\gtil,\Gtil_n(\Oo))$--action on $V_k^n$.}

\subsubsection{Proof.} The action of the $L_m^S$ given by the
Segal-Sugawara operators is well-defined on $V_k^n$ because $V^k_n$ is
a smooth $\ghat$--module. Next, we claim $L_m^S\cdot\vac_n=0$ for
$m\geq 2n-1$. For (precisely)
such $m$, for each term $:J^a_{m-i}J_{a,i}:$ at least one of the two
factors lies in $\gtil_n(\Oo)$, either immediately annihilating
$\vac_n$ or first passing through the other factor to annihilate
$\vac_n$, leaving a commutator of degree $m$ which also annihilates
$\vac_n$.

Observe that $\g_n(\Oo)$ acts locally nilpotently on $V_k^n$,
$Der_{2n}\Oo$ acts locally nilpotently on $\ghat/\g_n(\Oo)$ and
$\Der_{2n}\Oo$ annihilates $\vac_n$. This shows that
$\Der_{2n}\Oo$ acts locally nilpotently on $V_k^n$. It follows that
this action may be exponentiated to the pro-unipotent group
$\Aut_{2n}(\Oo)$. The arguments for $V_{k,c}^n$ are identical.

The fact that the $Vir$--action on $V_k^n$ defined by the
Segal-Sugawara operators is compatible with the action of $\ghat$
follows from commutation relations\eqref{LnJ}.

\subsubsection{The Segal-Sugawara singular vector.} \label{sug singular}

We now have two Virasoro actions on the representations $V_{k,c}^n$:
one given by the operators $L_m$ and one given by the operators
$L_m^S$. Moreover, both actions have the same commutation relations
with the Kac-Moody generators $J^a_n$.  Their difference now defines a
third Virasoro action, which has the crucial feature that it commutes
with $\ghat$. Indeed, setting $\Sug_m=L_m-L_m^S$, we have
\begin{eqnarray*}
[\Sug_l,\Sug_m]&=& [L_l,L_m]+ ([L_l^S,L_m^S]-[L_l,L_m^S]-
[L_l^S,L_m])\\ &=& [L_l,L_m]-[L_l^S,L_m^S]\\ &=&
(l-m)\Sug_{l-m}+\frac{c_k}{12}(l^3-l)\delta_{l,-m}
\end{eqnarray*}
and
$$[\Sug_m,J^a_l]=0,$$
where we have introduced the notation
\begin{equation}\label{c_k}
c_k=c-\frac{k\on{dim}\g}{k+\hv}
\end{equation}
for the central charge of the $\Sug_m$.

These operators are also defined from the action of a vertex operator.
Define the Segal-Sugawara singular vector
$\Sug\in V_{k,c}$ as the difference
$\Sug=\omega_V-\omega_S$, for $k\neq -h^{\vee}$.
The corresponding field $$\Sug(z)=\sum_{m\in\Z} \Sug_m z^{-m-2}$$
generates the action of the $\Sug_m$'s on $V_{k,c}^n$.
The crucial property of $\Sug$ is that it is a singular
vector for the Kac-Moody action, i.e., it is $\g(\Oo)$--invariant:
\begin{equation*}
J^a_n\cdot\Sug=0,\hskip.2in n\geq 0.
\end{equation*}
In what follows we consider $\Vck$ as a $\gtil$--module, where
$\ghat\subset \gtil$ acts by zero. By \propref{quadratic level},
$V_k^n$ is also a $\gtil$--module for $k \neq -h^\vee$. Therefore
their tensor product is a $\gtil$--module.

\subsubsection{\bf Proposition.} \label{Sugawara embedding}
{\em Let $k,c\in\C$ with $k\neq -\hv$.
\begin{enumerate}
\item The action of the $\Sug_m$ on $V_{k,c}^n$ defines an
embedding $\Sug^n_{k,c}: Vir_{c_k}^{2n} \longrightarrow V_{k,c}^n$
of $Vir$--modules.

\item $\Sug^n_{k,c}$ is a homomorphism of $\Gtil_n(\Oo)$--modules with
respect to the standard action of $\Gtil_n(\Oo)$ on $V_{k,c}^n$
and the trivial action of $G_n(\Oo)$ on $Vir_{c_k}^{2n}$.

\item There is an isomorphism of $(\gtil,\Gtil_n(\Oo))$--modules
\begin{equation*}
\sug_{k,c}:V_k^n \und{\C}\ot \Vck\longrightarrow V_{k,c}^n
\end{equation*}
such that $\sug_{k,c}(\vac\ot v)=\Sug^n_{k,c}(v)$.
\end{enumerate}
}

\subsubsection{Proof.}
The morphism $\Sug_{k,c}^n$ is uniquely determined by the requirement
that it intertwine the Virasoro action on $\Vck$ with the action of
the Virasoro algebra generated by the $\Sug_m$'s on $V_{k,c}^n$, which
annihilate the vacuum vector for $m\geq 2n-1$. (In particular, for
$n=0$, $\omega_V$ is sent to the Segal-Sugawara singular vector
$\Sug=\Sug_{-2}\cdot \vac_n$.) To show that this is an embedding, note
that for $m<2n-1$ the operators $L_m$ act freely on $V_{k,c}^n$, hence
so do the $\Sug_m$'s, which have the same leading term with respect to
the filtration on $V_{k,c}^n$ by the order in the $L_m$ operators.


The Segal-Sugawara operators commute with the action of $\ghat$. It
follows that the subspace $\Sug_{k,c}(Vir_{c_k}^{2n})$ of
$V_{k,c}^n$ generated by the action of the $\Sug_m$'s on $\vac_n$ is
$\g_n(\Oo)$--invariant, and hence $G(\Oo)$--invariant. Therefore any
vector $v\in\Sug_{k,c}^n(Vir_{c_k}^{2n}) \subset V_{k,c}^n$
determines a unique $\ghat$--homomorphism $v_*:V_k^n\to V_{k,c}^n$
with $v_*(\vac_n)=v$.  Hence we have a natural embedding of
$\Sug_{k,c}^n(Vir_{c_k}^{2n})$ in
$\on{Hom}_{\ghat}(V_k^n,V_{k,c}^n)$, and therefore a
$\ghat$--homomorphism
\begin{equation}\label{singular homomorphism}
\sug_{k,c}^n:V_k^n\ot \Vck \longrightarrow V_{k,c}^n.
\end{equation}
The fact that this is a homomorphism of $Vir$--modules, and therefore
of $\gtil$--modules, of this map is immediate from the formula
$\Sug_n=L_n-L_n^S$ for the Virasoro action, where $L_n^S$ and $L_n$
denote the Virasoro actions on $V_k^n$ and $V_{k,c}^n$,
respectively. In particular, we see that the map $\Sug_{k,c}^n$,
identified with the inclusion $\Vck\to \vac\ot\Vck$, followed by
$\sug_{k,c}^n$, is a homomorphism of $\Gtil_n(\Oo)$--modules, since
the actions of $\Gtil_n(\Oo)$ fixes the vector $\vac\subset V_k^n$.

We claim that the map $\sug_{k,c}^n$ in \eqref{singular homomorphism}
is an isomorphism.  By the Poincar\'e--Birkhoff--Witt theorem,
$V_{k,c}^n$ has a basis of monomials of the form
$$J^{a_1}_{m_1}\cdots J^{a_l}_{m_k}L_{n_1}\cdots L_{n_l} \vac_n
\hskip.3in(m_i<n,\; n_j<2n-1),$$
where we choose an ordering on the Kac-Moody and Virasoro generators.
It follows that the same holds with the $L_{n_j}$ replaced by
$\Sug_{m_j}$. The map $\sug_{k,c}^n$ acts as follows:
$$ (J^{a_1}_{m_1}\cdots J^{a_l}_{m_k}\vac)\ot (L_{n_1}\cdots
L_{n_l}\vac) \mapsto J^{a_1}_{m_1}\cdots J^{a_k}_{m_k}
\Sug_{n_1}\cdots\Sug_{n_l}\vac_n.$$ Hence it is indeed an isomorphism
as claimed.

\subsection{Limits of Segal-Sugawara.}

We would like to describe the behavior of the homomorphisms
$\Sug_{k,c}^n$, or equivalently, of the vector $\Sug\in V_{k,c}^n$
which generates them, as we vary the parameters $k,c$. In order to
describe the different limits, it is convenient to introduce the
parameters $\la,\mu$ and consider the full four--parameter family
$V_{k,c,\la,\mu}$ of vertex algebras. The vector
$\Sug=\omega_V-\omega_S$ is a well-defined element of
$V_{k,c,\la,\mu}$ for $k\in\C\sm -\hv$, $c\in\C$ and $\la,\mu\in\C\sm
0$.  It has a first order pole when $k=-h^{\vee}$, since the
Segal-Sugawara central charge $c_k=c-\dfrac{k\on{dim}\g}{k+\hv}$ does.
It also has a second order pole along $\la=0$ since it contains a term
quadratic in the $J$ generators, a first order pole along $\mu=0$ and
$c=\infty$ since it is first order in the $L$ generators, and a first
order pole along $k=\infty$ since it is quadratic in the $J$'s but
divided by $k+h^{\vee}$. Thus in these limits it is necessary to
normalize $\Sug$ by its leading term to obtain well-defined, non-zero
limits of the map $\Sug_{k,c}$.

\subsubsection{Critical level.}\label{critical level}

We would like to specialize the vector $\Sug$ and morphism
$\Sug_{k,c}$ to the critical level $k=-\hv$.  Introduce the rescaled
operators $$\ol{\Sug}_m=(k+\hv)\Sug_m=(k+\hv)(L_m-L_m^S),$$ which are
generated by the vertex operator $\ol{\Sug}(z)$ associated to
$$\ol{\Sug}=(k+\hv)\Sug=(k+\hv)L_{-2}-S_{-2}.$$ Thus when $k=-\hv$,
the vector $\ol{\Sug}$ is well-defined and equal to $-S_{-2}$.  The
operators $\ol{\Sug}_m$ satisfy
\begin{equation}\label{critical Sug commutators}
[\ol{\Sug}_l,\ol{\Sug}_m]= (k+\hv)((l-m)\ol{\Sug}_{l+m}
+\frac{(k+\hv)(k\on{dim}\g)}{12}(l^3-l)\delta_{l,-m}).
\end{equation}
Let us introduce the notation
$$\mu_\g=\hv\on{dim}\g.$$ We see that as $k$ approaches the critical
level $-\hv$, as the central charge of the Virasoro action of the
$\Sug_m$ becomes infinite, the renormalized operators $\ol{\Sug}_m$
become commuting elements, and moreover satisfy the Poisson relations
of the classical Virasoro algebra
$$U_{\infty,\mu_\g}Vir=\on{Sym}_{\mu_\g}Vir.$$

\subsubsection{\bf Proposition.} \label{critical Sugawara}
{\em The action of the operators $\ol{\Sug}_m$ defines a
$G(\Oo)$--invariant homomorphism
$\ol{\Sug}_{-\hv,c}^n:\ol{Vir}_{\mu_\g}^{2n} \longrightarrow
V_{-\hv,c}^n$ of $\on{Sym}_{\mu_\g}Vir$--modules. For $n=0$,
$\ol{\Sug}_{-\hv,c}$ defines a homomorphism of vertex Poisson algebras
$\ol{Vir}_{\mu_\g} \rightarrow {\mc Z}(V_{-\hv,c}(\g))$ to the center
of the Virasoro--Kac-Moody vertex algebra.}

\subsubsection{Proof.}
The morphism $\ol{\Sug}_{-\hv,c}^n$ is defined by the universal
mapping property of the vacuum module $\ol{Vir}_{\mu_\g}^n$ of
$\on{Sym}_{\mu_\g}Vir$.  The centrality of $\ol{\Sug}$ (and hence the
map) for $k=-\hv$ follows from the fact that the commutators of
$\ol{\Sug}_n=(k+\hv)\Sug_n$ (for $k\neq -\hv$) with the $L_n$ and
$J^a_n$ are divisible by $k+\hv$.  That this is a morphism of vertex
Poisson algebras is immediate from the commutation relation
\eqref{critical Sug commutators}.

\subsubsection{Infinite limit.}    \label{infinite limit}
Now we would like to study the ``generic'' classical limit of the
Segal-Sugawara construction in $\ol{V}_{\la,\mu}$. In order to do so
we approach the plane $c=k=\infty$ along the direction
$\dfrac{\la}{k}=\dfrac{\mu}{c}$, which is the direction we used to
define the vertex Poisson structure on $\ol{V}_{\la,\mu}(\g)$.  The
Segal-Sugawara operators are rescaled as follows:
$\ol{\Sug}_m=\frac{\la^2}{k}\Sug_m.$ These are the Fourier modes of
the vertex operator $\ol{\Sug}(z)$ associated to the vector
$$\ol{\Sug}=\frac{\la^2}{k}\Sug\in V_{k,c,\la,\mu},$$
which is regular for $k\neq 0,-\hv$ and $c,\la,\mu$ arbitrary.
In terms of the regular elements $\ol{L}$ and $\ol{J}$, we have
$$\ol{\Sug}=\la\ol{L_{-2}}-\frac{1}{2(k+h^{\vee})} \sum_a
\ol{J}^a_{-1} \ol{J}_{a,-1}+\cdots .$$ Thus when $\la=0$ the linear
term drops out and we recover the symbol,
$-\dfrac{1}{2(k+h^{\vee})}\ol{S}_{-2}$.  The commutation relations for
the $\ol{\Sug}_m$ are as follows:
\begin{eqnarray}\label{classical Sug commutators}
[\ol{\Sug}_l,\ol{\Sug}_m] &=&\frac{\la^4}{k^2}[\Sug_l,\Sug_m]\\
&=&\la\frac{\la}{k}((l-m)\ol{\Sug}_{l-m} +
\frac{\la\mu}{12}(l^3-l)\delta_{l,-m} +\cdots).
\end{eqnarray}
We have used the relation $\mu=\dfrac{\la c}{k}$, and that
$$\frac{\la^2}{k}\frac{k\on{dim}\g}{k+h^{\vee}} = \frac{\la^2
\on{dim}\g}{k+h^{\vee}}$$ vanishes in the limit
$k\to\infty$. Therefore in this limit the $\ol{\Sug}_m$'s satisfy the
relations of the Virasoro Poisson algebra $\Sym_{\la\mu}Vir$, with the
bracket rescaled by $\la$.  We will refer to a morphism which is a
homomorphism after rescaling by $\la$ as a $\la$--{\it
homomorphism}. We therefore obtain the following analogue of
\propref{critical Sugawara}:

\subsubsection{\bf Proposition.} \label{classical Sugawara}
{\em The action of the Segal-Sugawara operators $\ol{\Sug}_m$ defines a
$\la$--homo\-morphism $\ol{\Sug}_{\la,\mu}^n:\ol{Vir}_{\la\mu}^{2n}
\longrightarrow \ol{V}_{\la,\mu}^n$ of $\Sym_{\la\mu}Vir$--modules.
The image is $G(\Oo)$--invariant, and $\ol{\Sug}_{\la,\mu}$ is a
$\la$--homomorphism of vertex Poisson algebras.}

\section{Moduli Spaces.}\label{moduli spaces}





In this section we describe the spaces on which representations of the
Virasoro and Kac-Moody algebras and their semi-direct product
localize. These are the moduli spaces of curves, of bundles on a fixed
curve and bundles on varying curves, respectively. We will consider
these localization functors following the general formalism outlined
in \secref{patterns}. Note that the above moduli spaces are not
algebraic varieties, but algebraic {\em stacks}. However, as explained
in \cite{Jantzen,Hecke}, the localization formalism is applicable to
them because they are ``good'' stacks, i.e., the dimensions of their
cotangent stacks are equal to the twice their respective dimensions.

\subsection{Moduli of bundles.} \label{moduli of bundles}

Let $X$ be a smooth projective curve over $\C$. Denote by $\Bun_G(X)$
the moduli stack of principal $G$--bundles on $X$, and $\Pf$ the
tautological $G$--bundle on the product $X\times \Bun_G(X)$. (Its
restriction to $X\times\{\mc P\}$, for ${\mc P}\in \Bun_G(X)$, is
identified with $\mc P$.)

Given $x \in X$, we denote by $\Bun_G(X,x,n)$ the moduli stack of
$G$--bundles with an $n$th order jet of trivialization at $x$, and by
$\wh{\Bun}_G(X,x)$ the moduli stack of $G$--bundles with
trivializations on the formal neighborhood of $x$ (the latter moduli
stack is in fact a {\em scheme} of infinite type). For now we fix a
formal coordinate $t$ at $x$, so that the complete local ring $\Oo_x$
is identified with $\Oo = \C[[t]]$. Later we will vary this coordinate
by the $\AutO$--action. The group scheme $G(\Oo)$ acts on
$\wh{\Bun}_G(X,x)$ by changing the formal trivialization, making
$\wh{\Bun}_G(X,x)\to \Bun_G(X)$ into a $G(\Oo)$--torsor.  More
generally $\wh{\Bun}_G(X,x)\to \Bun_G(X,x,n)$ is a $\Gn$--torsor.

\subsubsection{\bf Theorem.}\label{Kac-Moody uniformization}
(Kac-Moody Uniformization.)
{\em
\begin{enumerate}
\item The $\Gn$ action on the moduli space $\wh{\Bun}_G(X,x)$ extends
to a formally transitive action of the Harish-Chandra pair
$(\g(\K),\Gn)$.

\item (\cite{BL,BL2,DSi}) For $G$ semisimple, the action of the
ind-group $G(\K)$ on $\wh{\Bun}_{G}(X,x)$ is transitive, and there are
isomorphisms
$$\wh{\Bun}_{G}(X,x) \simeq G(\K)_{\on{out}} \backslash G(\K), \quad
\quad \Bun_G(X) \simeq G(\K)_{\on{out}} \backslash G(\K)/G(\Oo).$$
\end{enumerate}
}

\subsubsection{Line bundles on $\Bun_G(X)$.}
We refer to \cite{Sor} for a detailed discussion of the line bundles
on $\Bun_G(X)$. They are classified by integral invariant forms on
$\g$, which also label the Kac-Moody central extensions of $LG$.  The
action of $\ghat$ on $\wh{\Bun}_G(X,x)$ lifts with level one to the
line bundle $\Cc$ given by the corresponding invariant form.  This
line bundle may be defined by using \thmref{Kac-Moody uniformization}
from the action of the Kac-Moody {\em group} $\wh{G(\K)}$ (the central
extension splits over $G(X\sm x)$ and hence gives rise to a line
bundle on $G(X\sm x)\bs \wh{G(\K)}=\wh{\Bun}_G(X,x)$, which descends to
$\Bun_G(X)$).

For example, if $G=SL_n$, the line bundle $\Cc$ may be identified with
the determinant of the cohomology of the universal vector bundle ${\mf
E}=\Pf \und{SL_n}\times\C^n$ over $X\times \Bun_{SL_n}(X),$
$\on{det}R\pi_{2*}{\mf E}$ (where $\pi_2:X\times \Bun_G(X)\to
\Bun_G(X)$). This identification however is not canonical (it is not
valid for bundles over varying curves, see \secref{moduli of curves
and bundles}). More generally, for any simple algebraic group $G$
powers of $\Cc$ can be defined as determinant line bundles associated
to representations of $G$.

\subsubsection{Localization.}

For every $n$ the triple
$$(\Mhat,\M,\Ll)=(\wh{\Bun}_G(X,x), \Bun_G(X,x,n),\Cc)$$ defined above
carries a transitive Harish--Chandra action of $(\ghat,\Gn)$ as in
\secref{localization}. Therefore we have localization functors from
$(\ghat,\Gn)$--modules to twisted $\D$--modules on $\Bun_G(X,x,n)$ (we
denote these functors by $\Delta$ as before). In particular, according
to \propref{twists are algebras}, (2), for the vacuum module $V_k^n$,
the sheaf $\Delta(V_k^n)$ on $\Bun_G(X,x,n)$ is just the corresponding
sheaf of twisted differential operators, which we denote uniformly by
$\D_k$. Furthermore, the twist $\V_k^n = \V_{\ghat,G_n(\Oo)}$ is a
sheaf of algebras on $\Bun_G(X,x,n)$, and we have a surjective
homomorphism $\V_k^n \to \D_k$.

The corresponding classical vacuum representations $\ol{V}_k^n$
localize to give the Poisson sheaves of functions on the twisted
cotangent bundles $T^*_{\la}\Bun_G(X,x,n)=\Conn\Cc^{\la}$
corresponding to $\Cc$. Recall that the cotangent space
$T^*_{\Pc}\Bun_G(X)$ at a bundle $\Pc$ is the space of $\g_{\mc
P}$--valued differentials $\on{H}^0(X,\g_{\mc
P}\otimes\Omega)$. Therefore we obtain

\subsubsection{\bf Proposition. (\cite{Fa,BS})} \label{conn as TCB}
{\em The twisted cotangent bundle $T^*_1 \Bun_G(X,x,n)$ is canonically
identified, as a torsor over $T^*\Bun_G(X,x,n)$, with the moduli stack
\linebreak $\Conn_G(X,x,n)$ of bundles with connections with a pole of
order at most $n$ at $x$. }

\subsection{Moduli of curves.}\label{moduli of curves}
Let $\M_g$ denote the moduli stack of smooth projective curves of
genus $g$, and $\pi:{\Xf_g}\to\M_g$ the universal curve. The stack
$\Xf_g$ is identified with the moduli stack $\Mgo$ of {\em pointed}
genus $g$ curves. More generally, we denote by $\M_{g,1,n}$ the moduli
stack of pointed curves with an $n$th order jet of coordinate at the
marked point, and $\Mghat=\M_{g,1,\infty}$ the moduli scheme of curves
with a marked point and formal coordinate (i.e. classifying triples
$(X,x,z)$, where $(X,x)\in\Mgo$ and $z$ is a formal coordinate on $X$
at $x$). The group scheme $\AutO$ acts on $\Mghat$ by changing the
coordinate $z$, making $\Mghat$ into an $\AutO$--torsor over $\Mgo$
and an $\An$--torsor over $\M_{g,1,n}$.

For any family of curves $\pi:X\to S$, we have the Hodge line bundle $\Hh$
on $S$, defined by
$$\Hh=\on{det}R\pi_*\omega_{X/S},$$ the determinant of the cohomology
of the canonical line bundle of $X$ over $S$. By abuse of notation we
will denote by $\Hh$ the Hodge line bundle of an arbitrary family of
curves. Over the moduli stack $\M$, Mumford has shown that (for $g>1$)
$\Hh$ generates the Picard group $\on{Pic}\M\cong\Z\cdot\Hh$.

\subsubsection{\bf Theorem.}\label{Virasoro uniformization} \cite{BS,TUY,ADKP,Kont}
(Virasoro Uniformization.)  {\em The $\AutO$ action on the moduli
space $\Mghat$ of pointed, coordinatized curves extends to a formally
transitive action of the Harish-Chandra pair $(Vir,\AutO)$ (of level
$0$).  The action of $Vir$ lifts to an action with central charge $-2$
on the line bundle $\Hh$.}

\subsubsection{Localization.}\label{localization for Vir}

The pattern of localization from \secref{localization} applies
directly to the $\Aut_n\Oo$--bundle $\Mghat\to \M_{g,1,n}$ and the
Harish--Chandra pairs $(Vir,\Aut_n\Oo)$. Therefore we have
localization functors from $(Vir,\Aut_n\Oo)$--modules to twisted
$\D$--modules on $\M_{g,1,n}$. In particular, the localization of the
vacuum module $Vir_c^n$ gives us the sheaf of twisted differential
operators on $\M_{g,1,n}$, denoted by $\D_c$, and their classical
versions give us Poisson algebras of functions on twisted cotangent
bundles. Furthermore, the twist $\Vir_c^n = \V_{Vir,\Aut_n(\Oo)}$ is a
sheaf of algebras on $\M_{g,1,n}$, and we have a surjective
homomorphism $\Vir_c^n \to \D_c$.

A similar picture holds for the fibration $\Mghat\to \M_g$ and the
pair $(Vir,\Der\Oo)$, except that $\Der\Oo$ does not integrate to a
group (only to an ind-group, see \cite{book}, Ch. 5), and $\Mghat\to
\M_g$ is not a principal bundle, so that the definition of
localization does not carry over directly. Nevertheless, we obtain the
desired description of differential operators on $\M_g$ by first
localizing $Vir_c$ on $\Mgo$ (or $\M_{g,1,n}$), using the fact that
the corresponding $\D$--modules descend along the projection
$\pi:\Mgo\to\M_g$:

\subsubsection{\bf Proposition.}\label{pullback of diffops}
{\em The localization of the Virasoro module $Vir_c$ on $\Mgo$ is
isomorphic to the pullback $\pi^*\D_c$ of the sheaf of twisted
differential operators on $\M_g$.}

\subsubsection{Proof.}
By \lemref{change of vacuum}, the localization $\Delta(Vir_c)$ is the
quotient of $\Delta(Vir_c^0) = \D_c$ by the action of the partial
vacuum representation $\Delta(V_{\Der \Oo,\Der_0\Oo})$.  However, the
latter is readily identified as the sheaf of {\em relative}
differential operators, $\Delta(V_{\Der
\Oo,\Der_0\Oo})=\D_{c/\M_g}$. Indeed the action of $\Der\Oo$ on
$\wh{\M}_g$ is free and generates the relative vector fields on the
universal curve -- the $\wh{\M}_g$ twist
$(\Der\Oo/\Der_0\Oo)_{\wh{\M}_g}$ is precisely the relative tangent
sheaf of $\Mgo$ over $\M_g$.  But the quotient of the sheaf of
differential operators by the ideal generated by vertical vector
fields is the pullback of the sheaf of differential operators
downstairs, whence the proposition.

\subsubsection{}
Recall (see, e.g., \cite{book}) that the space of projective
structures on $X$ is a torsor over the quadratic differentials
$\on{H}^0(X,\Omega^{\otimes 2})$, which is the cotangent fiber of
$\M_g$ at $x$. The Virasoro uniformization of $\M_g$, together with
the canonical identification of projective structures on the punctured
disc with a hyperplane in $Vir^*$, give us the following:

\subsubsection{\bf Corollary.} (\cite{BS})\label{proj as TCB}
{\em There is a canonical identification of the twisted cotangent
bundles $T^*_{\la}\M_g=\Proj_g$ (for $\la=12$).}

\subsubsection{}Similarly the twisted cotangents to the moduli $\M_{g,1,n}$ of curves
with marked points and level structures are identified with the moduli
$\Proj_{g,1,n}^{\la}$ of $\la$--projective structures with poles at the
corresponding points.

\subsection{Moduli of curves and bundles.} \label{moduli of curves and
bundles}

The discussions of the moduli spaces $\M_g$ and $\Bun_G(X)$ above may
be generalized to the situation where we vary both the curve and the
bundle on it. Let $\Bun_{G,g}$ denote the moduli stack of pairs
$(X,\Pc)$, where $X$ is a smooth projective curve of genus $g$, and
$\Pc$ is a principal $G$--bundle on $X$. Thus we have a projection
$$\Pi:\Bun_{G,g} \longrightarrow\M_g$$ with fiber over $X$ being the
moduli stack $\Bun_G(X)$. Let $\Bun_{G,g,1}$ denote the moduli stack
classifying $G$--bundles on pointed curves, i.e., the pullback of the
universal curve to $\Bun_{G,g}$:
$$\Bun_{G,g,1}=\Mgo \und{\M_g}\times \Bun_G \stackrel{\pi_G}{\lrarrow}
\Bun_{G,g}$$ We denote by $\Pi^0:\Bun_{G,g,1}\to \Mgo$ the map
forgetting the bundles.  We let $\Bun_{G,g,1,n}$ denote the moduli
stack of quintuples $(X,x,\Pc,z,\tau)$ consisting of a $G$--bundle
$\Pc$ with an $n$th order jet of trivialization $\tau$ on a pointed
curve $(X,x)$ with $2n$th order jet of coordinate $z$. We let
$\Pi^n:\Bun_{G,g,1,n}\to \M_{g,1,2n}$ denote the map forgetting
$(\Pc,\tau)$. For $n=\infty$ we obtain the moduli scheme
$\wh{\Bun}_{G,g}$ of bundles with formal trivialization on pointed
curves with formal coordinates.

There is a natural two--parameter family of line bundles on
$\Bun_{G,g}$. Namely, there is a Hodge bundle $\Hh$ associated to the
family of curves $\pi_G:\Bun_{G,g,1}\to \Bun_{G,g}$ (which is the
$\Pi$ pullback of the Hodge bundle on $\M_g$).  The universal
principal $G$--bundle $\Pf$ lives on this universal curve
$\Bun_{G,g,1}$, so we may also consider the line bundle $\Cc$
associated to the principal bundle $\Pf$.  The bundle $\Cc$ is trivial
on the section $\on{triv}:\M_g\to \Bun_{G.g}$ sending a curve to the
trivial $G$--bundle. Let $\Ll_{k,c}=\Cc^{\ot k}\ot \Hh^{\ot c}$. For
simply-connected $G$ this assignment gives an identification
$\on{Pic}(\Bun_{G,g})\cong \Z\oplus\Z$ (see \cite{Laszlo}).

For $G=SL_n$ we may consider the determinant bundle
$\on{det}R\pi_{2*}{\mf E}$ as before, where $\pi_2$ is the projection
from the universal curve and ${\mf E}=\Pf\times_{SL_n}\C^n$ is the
universal vector bundle. The determinant of the cohomology of the
trivial rank $n$ bundle gives the $n$th power of the Hodge line
bundle, while for fixed curve the determinant bundle may be identified
with $\Cc$. Thus we have the Riemann-Roch identification
$$\on{det}R\pi_{2*}{\mf E}\cong\Cc\ot\Hh^{\ot n}=\Ll_{1,n}.$$ (In this
case the determinant and Hodge bundles also span the Picard group.) In
general we have determinant line bundles $\Ll_{\rho}$ for any
representation $\rho:G\to SL_n$.

\subsubsection{Extended connections}

For each $\la,\mu\in\C$, we have the corresponding twisted cotangent
bundle $T^*_{\la,\mu}\Bun_{G,g}$, which is $\Conn(\Ll_{\la,\mu})$ when
$\la,\mu$ are integral. Following \cite{szego,thetas}, we will refer to
points of $\ExConn^{\la,\mu}_{G,g}=T^*_{\la,\mu}\Bun_{G,g}$ as
$(\la,\mu)$--{\em extended connections}.
Using the short exact
sequence of cotangent bundles
$$0\to \Pi^*T^*\M_g\to T^*\Bun_{G,g}\to T^*_{/\M}\Bun_{G,g}\to 0$$ and
the description (see \propref{conn as TCB}) of the relative twisted
cotangent bundle of $Bun_{G,g}$, we obtain an affine projection
\begin{equation}\label{exconn to conn}
\ExConn^{\la,\mu}_{G,g}\to\Conn^{\la}_{G,g},
\end{equation} with fibers affine
spaces over quadratic differentials.  Here $\Conn_{G,g}$ denotes the
moduli space of curves equipped with $G$--bundles and
$\la$--connections, which is identified as the {\em relative} twisted
cotangent bundle
$$\Conn^{\la}_{G,g}\cong \Conn_{/\M}\Ll_{\la,\mu}$$
for any $\mu\in\C$. (Note that $\Ll_{0,\mu}=\Pi^*\Hh_\mu$ has a
canonical connection relative to $\M$, so that
$\Conn_{/\M}\Ll_{\la,\mu} \cong \Conn_{/\M}\Ll_{\la,0}$ for any $\mu$.)

\subsubsection{Remark: kernel functions.}
See also \cite{szego} where a concrete description is given of the
identification of projective structures, connections and extended
connections with the twisted cotangent bundles of the moduli of curves
and bundles using kernel functions along the diagonal (in particular
the Szeg\"o kernel), in the spirit of \cite{BS}.

\medskip

\subsubsection{}
The group schemes $\Gtil_n(\Oo)$ act on $\wh{\Bun}_{G,g}$, changing the
coordinate $z$ and trivialization $\tau$. This action makes
$\wh{\Bun}_{G,g}$ into a $\Gtil_n(\Oo)$--torsor over $\Bun_{G,g,1,n}$.

\subsubsection{\bf Theorem.}\label{Virasoro-Kac-Moody uniformization}
(Virasoro--Kac-Moody Uniformization.)  {\em The $\Gtil_n(\Oo)$--action
on the moduli stack $\wh{\Bun}_{G,g}$ extends to a formally transitive
action of the Harish-Chandra pair $(\gtil, \Gtil_n(\Oo))$ (of level
and central charge $0$). The action of $\gtil$ lifts to an action
with level $k$ and central charge $c$ on the line bundle $\Ll_{k,c}$.}

\subsubsection{Localization.} The pattern of localization from
\secref{localization} again applies directly to $\Bun_{G,g,1,n}$ and
the Harish--Chandra pairs $(\gtil,\Gtil_n(\Oo))$. Therefore we have
localization functors from $(\gtil,\Gtil_n(\Oo))$--modules to twisted
$\D$--modules on $\Bun_{G,g,1,n}$. In particular, the localization of
the vacuum module $V_{k,c}^n$ gives us the sheaf of twisted
differential operators on $\Bun_{G,g,1,n}$, denoted by $\D_{k,c}$. As
in \secref{localization for Vir}, we would like to show that the
corresponding sheaves descend along the projection
$\pi_G:\Bun_{G,g,1}\to \Bun_{G,g}$ to the moduli of curves and bundles
itself.  The argument of \propref{pullback of diffops} carries over
directly to this setting, and we obtain

\subsubsection{\bf Proposition.}\label{pullback of diffops for bundles}
{\em The localization of the Virasoro--Kac-Moody module $V_{k,c}$ on
\linebreak $\Bun_{G,g,1}$ is isomorphic to the pullback
$\pi_G^*\D_{k,c}$ of the sheaf of twisted differential operators on
$\Bun_{G,g}$.}

\bigskip

Note that we now have a {\em two}--parameter family of line bundles
$\Ll_{k,c}$, and hence the pattern of deformations of \secref{limits}
can be ``doubled'', to match up with the picture of the
Virasoro--Kac-Moody vertex algebra $V_{k,c}$. Thus, we introduce
deformation parameters $\la,\mu$ coupled to the level and charge
$k,c$. This defines a four--parameter family of algebras
$\D_{\la,\mu}(\Ll_{k,c})$, to which the analogous quasi-classical
localization statements apply.

\section{The Segal-Sugawara Homomorphism.}\label{Sugawara section}

In this section we apply the techniques of Lie algebra localization
and vertex algebra conformal blocks to the Segal-Sugawara construction
from \secref{Vir-KM algebras}. We interpret the result as a
homomorphism between sheaves of twisted differential operators on the
moduli stacks introduced in the previous section, and as heat
operators on spaces of nonabelian theta functions. Various classical
limits of this construction will be considered in the next section.

\subsection{Homomorphisms between sheaves of differential operators.}

Let $k,c\in\C$ with $k\neq -\hv$. Let $\V_{k,c}^n =
\V_{\gtil,\Gtil_n(\Oo)}$ be the sheaf on $\Bun_{G,g,1,n}$ obtained as
the twist of the vacuum module $V_{k,c}^n$ over $(\gtil,\Gtil_n(\Oo))$
following the construction of \secref{localizing vacuum}. By
\propref{twists are algebras}, this is a sheaf of algebras, equipped
with a surjective homomorphism to the sheaf of $(k,c)$--twisted
differential operators $\Delta(V_{k,c}^n)=\D_{k,c}$.

According to \corref{global diff}, the subspace of invariants
$(V_{k,c}^n)^{\Gtil_n(\Oo)}$ give rise to {\em global} differential
operators on $\Bun_{G,g,1,n}$. Unfortunately, for $k \neq -h^\vee$ the
space $(V_{k,c}^n)^{\Gtil_n(\Oo)}$ is one-dimensional, spanned by the
vacuum vector. However, the Segal-Sugawara construction provides us
with a large space of invariants for the smaller group $\Gn \subset
\Gtil_n(\Oo)$. Namely, by \propref{Sugawara embedding}, the
$\Gn$--invariants contain a copy of the Virasoro vacuum module
$\Sug^n: Vir_{c_k}^{2n}\to (V_{k,c}^n)^{\Gn}$. We will use this fact
to obtain a homomorphism $\D_{c_k} \to \Pi^n_* \D_{k,c}$ of
sheaves on $\M_{g,1,2n}$. The first step is the following assertion.

Consider $Vir_{c_k}^{2n}$ as a $\gtil$--module by letting $\ghat$
act by zero. Then the twist $\Vir_{c_k}^{2n} =
\V_{\gtil,\Gtil_n(\Oo)}$ becomes a sheaf of algebras on $\M_{g,1,2n}$.

\subsubsection{\bf Proposition.}\label{homomorphism on twists}
{\em There is a homomorphism of sheaves of algebras
on $\M_{g,1,2n}$,
$$\Sug^n_{k,c}:\Vir_{c_k}^{2n}\longrightarrow
\Pi_*^n\V_{k,c}^n.$$
}

\subsubsection{Proof.}
Note that the map $\Sug^n: Vir_{c_k}^{2n}\to (V_{k,c}^n)^{\Gn}$ is a
homomorphism of $\Gtil_n(\Oo)$--modules (where the subgroup $G_n(\Oo)$
of $\Gtil_n(\Oo)$ acts by zero). Hence its gives rise to a
homomorphism of the corresponding twists by the $\Gtil_n(\Oo)$--torsor
$\wh{\Bun}_{G,g} \to \Bun_{G,g,1,n}$,
\begin{multline}\label{algebra sequence}
\wh{\Bun}_{G,g} \und{\Bun_{G,g,1,n}}\times Vir_{c_k}^{2n}
\longrightarrow \wh{\Bun}_{G,g} \und{\Bun_{G,g,1,n}}\times
(V_{k,c}^n)^{\Gn} \longrightarrow \\ \longrightarrow \V_{k,c}^n =
\wh{\Bun}_{G,g} \und{\Bun_{G,g,1,n}}\times V_{k,c}^n.
\end{multline}
The first two sheaves are pullbacks from $\M_{g,1,2n}$. Indeed, since
the actions of $\Gtil_n(\Oo)$ on $Vir^{2n}_{c_k}$ and
$(\V_{k,c}^n)^{\Gn}$ factor through $\Atn$, their twists depend only
on the associated $\Atn$--torsor, which is nothing but the
$\Pi^n$--pullback of the $\Atn$--torsor $\pi_{2n}:
\Mghat\to\M_{g,1,2n}$. Therefore the maps \eqref{algebra sequence}
may be written as
$$
(\Pi^n)^*\Vir_{c_k}^{2n}\longrightarrow
(\Pi^n)^*(\V_{k,c}^n)^{\Gn} \longrightarrow \V_{k,c}^n.
$$

By adjunction, we obtain the following maps on $\M_{g,1,2n}$:
$$
\Vir_{c_k}^{2n}\longrightarrow (\V_{k,c}^n)^{\Gn} \longrightarrow
\Pi^n_* \V_{k,c}^n.
$$ We claim that these sheaves are algebras and these maps are algebra
homomorphisms. The third term is the pushforward of an algebra, hence
an algebra. The first term is the twist of the vacuum module for
$(Vir,\Atn)$, so its algebra structure comes from \propref{twists are
algebras}. Next, note that the Harish-Chandra pair
$(\gtil,\Gtil_n(\Oo))$ acts on $\Mghat$ through the quotient
$(Vir,\Aut_{2n} \Oo)$, so that the $\Atn$--twist of the
$\Gn$--invariants in $V_{k,c}^n$ may be rewritten in terms of the
semidirect product $\Gtil_n(\Oo)$:
$$(\V_{k,c}^n)^{\Gn} = \left[ \pi_{n*}(\Oo_{\Mghat}\otimes
V_{k,c}^n) \right]^{\Gtil_n(\Oo)}.$$ Thus as in \propref{twists are
algebras} the sheaves $(\V_{k,c}^n)^{\Gn}$ on $\M_{g,1,2n}$ or
$\Bun_{G,g,1,n}$ are naturally sheaves of algebras. This structure is
clearly compatible with that on $\Pi_*^n\V^n_{k,c}$.

Finally, the maps $Vir_{c_k}^{2n}\to V_{k,c}^n$ are induced from the
homomorphism $U(Vir_{c_k}) \to \wh{U}_{k,c}\gtil$ into the completion
of the enveloping algebra of $\gtil$. This homomorphism maps the
subalgebra $U\Der_{2n}\Oo$ to the left ideal generated by the Lie
subalgebra $\gtil_n(\Oo)$. Hence we obtain a homomorphism of the
corresponding vacuum modules, because $V_{k,c}^n$ can be defined as
the quotient of the completed algebra $\wh{U}_{k,c}\gtil$ by the left
ideal generated by the Lie subalgebra $\gtil_n(\Oo)$. It follows that
the map $\Vir_{c_k}^{2n} \longrightarrow (\V_{k,c}^n)^{\Gn}$ above
respects algebra structures.  More precisely, this map comes from the
above homomorphism of enveloping algebras by passing to the
normalizers of ideals on both sides, hence it remains a homomorphism.

\subsubsection{}
By \propref{homomorphism on twists}, we have a diagram of algebra
homomorphisms on $\M_{g,1,2n}$,
\begin{equation*}
\begin{array}{ccc}
\Vir_{c_k}^n&\stackrel{\Sug^n_{k,c}}\longrightarrow&
\Pi_*^n\V_{k,c}^n\\ \downarrow&&\downarrow\\ \D_{c_k}&&\Pi_*^n\D_{k,c}
\end{array}
\end{equation*}
We wish to show that the homomorphism $\Sug^n_{k,c}$ descends to the
sheaves of twisted differential operators. It will then automatically
be a homomorphism of {\em algebras} of differential operators.

\subsubsection{\bf Theorem.}\label{descent to coinvariants}
{\em The homomorphism $\Sug^n_{k,c}$ of descends to a
homomorphism of algebras
$$
\Sug^n_{k,c}:\D_{c_k}\longrightarrow\Pi_*^n\D_{k,c}.
$$
}

\subsubsection{Spaces of coinvariants.}    \label{reminder}

Equivalently, we need to show that up on $\Bun_{G,g,1,n}$, the
morphism $(\Pi^n)^* \Vir_{c_k}^n\to \V_{k,c}^n$ descends to
$(\Pi^n)^*\D_{c_k}\to\D_{k,c}$, since the morphism on $\M_{g,1,2n}$ is
obtained from the former by adjunction.

If we apply the Virasoro--Kac-Moody localization functor $\Delta$ on
$\Bun_{G,g,1,n}$ to the $(\gtil,\Gtil_n(\Oo))$--modules $Vir_{c_k}^{2n}$
and $V_{k,c}^n$, we obtain the desired sheaves $(\Pi^n)^*\D_{c_k}$ and
$\D_{k,c}$. However, the embedding $\Vck\to V_{k,c}$ of
\propref{Sugawara embedding} is {\em not} a homomorphism of
$\gtil$--modules: it intertwines the Virasoro action on $Vir_{c_k}^n$
with the action on $V_{k,c}$ of the Virasoro algebra generated by the
$\Sug_n$'s, not the $L_n$'s. Because of that, it is not immediately
clear that the map $\Vck\to V_{k,c}$ gives rise to a morphism of
sheaves $(\Pi^n)^*\D_{c_k}\to\D_{k,c}$. In order to prove that, we
must use \thmref{coinv} below to pass from Lie algebra coinvariants to
vertex algebra coinvariants.

Let us recall some results from \cite{book} on the spaces of (twisted)
coinvariants of vertex algebras.

Let $V$ be a conformal vertex algebra with a compatible
$\ghat$--structure (see \cite{book}, Section 6.1.3). This means that
$V$ carries an action of the Harish-Chandra pair $(\gtil,\Gtil(\Oo))$,
such that the action of the Lie algebra $\gtil$ is generated by
Fourier coefficients of vertex operators. In particular, our vacuum
modules $V_{k,c}$ and $Vir_c$ are examples of such vertex
algebras.

Given a vertex algebra of this type, we define its space of twisted
coinvariants as in \cite{book}, Section 8.5.3. Namely, let $R$ be a
local $\C$--algebra and $(X,x,\Pc)$ an $R$--point of $\Bun_{G,g}$,
i.e., a pointed curve $(X,x)$ and a $G$--bundle on $X$, all defined
over $\on{Spec} R$. Then $X$ carries a natural $\Gtil(\Oo)$--bundle
$\wh{\Pc}$, whose fiber over $y \in X$ consists of pairs $(z,t)$,
where $z$ is a formal coordinate at $y$ and $t$ is a trivialization of
$\Pc$ over the formal disc around $y$.

Set
$$
\V^{\Pc} = \wh{\Pc} \und{\Gtil(\Oo)}\times V.
$$
This vector bundle carries a flat connection, and we define the sheaf
$h(\V^{\Pc} \otimes \Omega)$ as the sheaf of zeroth de Rham cohomology
of $\V^{\Pc} \otimes \Omega$. The vertex algebra structure on $V$
makes this sheaf into a sheaf of Lie algebras. In particular, the Lie
algebra $U^{\Pc}(\V_x)$ of sections of $h(\V^{\Pc} \otimes \Omega)$
over the punctured disc $D_x^\times$ is isomorphic to the Lie algebra
$U(V)$ topologically spanned by the Fourier coefficients of all vertex
operators from $V$. It acts on $\V^{\Pc}_x$, the fiber of $\V^{\Pc}$
at $x$.

Let $U_{X\bs x}^\Pc(\V_x)$ be the image of the Lie algebra
$\Gamma(X\bs x,h(\V^{\Pc} \otimes \Omega))$ in $U^{\Pc}(\V_x) =
\Gamma(D_x^\times,h(\V^{\Pc} \otimes \Omega))$. The space
$H^{\Pc}(X,x,V)$ of twisted coinvariants of $V$ is by definition the
quotient of $\V^\Pc_x$ by the action of $U_{X\bs x}^\Pc(\V_x)$.

Let ${\mc A}^{\mc P}$ be the Atiyah algebroid of infinitesimal
symmetries of ${\mc P}$. We have the exact sequence
$$
0 \to \g^{\mc P} \to {\mc A}^{\mc P} \to \Theta_X \to 0,
$$ where $\g^{\mc P}$ is the sheaf of sections of the vector bundle
${\mc P} \und{G}\times \g$.

When $V = V_{k,c}$, the Lie algebra $U^\Pc(\V_{k,c,x})$ contains as a
Lie subalgebra a canonical central extension of the Lie algebra
$\Gamma(D_x^\times,{\mc A}^{\mc P})$ (it becomes isomorphic to $\gtil$
if we choose a formal coordinate at $x$ and a trivialization of ${\mc
P}|_{D_x}$). Also, the Lie algebra $U_{X\bs x}^\Pc(\V_{k,c,x})$
contains
$$
\gtil^\Pc_{\on{out}} = \Gamma(X \bs x,{\mc A}^{\mc P})
$$
as a Lie subalgebra (it is isomorphic to $\on{Vect}(X \bs x)
\ltimes \g \otimes \C[X\bs x]$).

Likewise, in the case when $V = Vir_c$, $U^\Pc(\Vir_{c,x})$ contains
the Virasoro algebra, and $U_{X\bs x}^\Pc(\Vir_{c,x})$ contains the
Lie algebra $\on{Vect}(X\bs x)$ of vector fields on $X\bs x$. The
homomorphism $Vir_{c_k} \to V_{k,c}$ of vertex algebra induces
injective homomorphisms of Lie algebras
\begin{align*}
U^\Pc(\Vir_{c_k,x}) & \hookrightarrow U^\Pc(\V_{k,c,x}) \\ U_{X\bs
x}^\Pc(\Vir_{c_k,x}) & \hookrightarrow U_{X\bs x}^\Pc(\V_{k,c,x})
\end{align*}
(since the action of $\ghat$ on $Vir_c$ is trivial, the Lie algebras
on the left do not depend on $\Pc$). Note that though the
Segal-Sugawara vertex operator is quadratic in the generating fields
of $V_{k,c}$, the elements of $U_{X\bs x}^\Pc(\Vir_{c_k,x})$ (and
other elements of $U_{X\bs x}^\Pc(\V_x)$) cannot be expressed in terms
of the Lie subalgebra $\gtil^\Pc_{\on{out}}$ of $U_{X\bs
x}^\Pc(\V_{k,c,x})$. Nevertheless, we have the following result which
is proved in \cite{book}, Theorem 8.3.3 (see also Remark 8.3.10).

\subsubsection{\bf Proposition.}    \label{coinv}
{\em For any smooth $\gtil^\Pc$--module $M$ (which is then
automatically a $U^\Pc(\V_{k,c,x})$--module), the space of
coinvariants of $M$ by the action of $U_{X\bs x}^\Pc(\V_{k,c,x})$ is
equal to the space of coinvariants of $M$ by the action of
$\gtil^\Pc_{\on{out}}$.}

\subsubsection{Proof of \thmref{descent to coinvariants}.}

Fix an $R$--point $(X,x,z,\Pc,t)$ of $\Bun_{G,g,1,n}$, i.e., a pointed
curve $(X,x)$, a $2n$--jet of coordinate $z$ at $x$, a $G$--bundle
$\Pc$ and an $n$--jet of trivialization $t$ of $\Pc$ at $x$, all
defined over the spectrum of some local $\C$--algebra $R$.  Let $v$ be
a vector in the fiber of the sheaf $\Vir^{2n}_{c_k}$ over the
$R$--point $(X,x,z)$ of $\M_{g,1,2n}$, which lies in the kernel of the
surjection $\Vir^{2n}_{c_k}|_{R} \to \D_{c_k}|_{R}$. In order to
prove the theorem, it is sufficient to show that the image of $v$ in
the fiber of the sheaf $\V^n_{k,c}$ over the $R$--point
$(X,x,z,\Pc,t)$ of $\Bun_{G,g,1,n}$ belongs to the kernel of the
surjection $\V^n_{k,c}|_{R} \to \D_{k,c}|_{R}$.

But according to \cite{book}, Lemmas 16.2.9 and 16.3.6, the kernel of
the map $\Vir^{2n}_{c_k}|_{R} \to \D_{c_k}|_{R}$ is spanned by all
vectors of the form $s \cdot A$, where $A \in \Vir^{2n}_{c_k}|_{R}$
and $s \in \on{Vect}(X\bs x)$ (so that $\D_{c_k}|_{R}$ is the space of
coinvariants of $\Vir^{2n}_{c_k}|_{R}$ with respect to
$\on{Vect}(X\bs x)$). Likewise, the kernel of the map $\V^n_{k,c}|_{R}
\to \D_{k,c}|_{R}$ is spanned by all vectors of the form $f \cdot B$,
where $B \in \V^n_{k,c}|_{R}$ and $f \in \gtil^{\Pc}_{\on{out}}$ (so
that $\D_{k,c}|_{R}$ is the space of coinvariants of
$\V^n_{k,c}|_{R}$ with respect to $\gtil^{\Pc}_{\on{out}}$).

So we need to show that the image of a vector of the above form $s
\cdot A$ in $\V^n_{k,c}|_{R}$ is in the image of the Lie algebra
$\gtil^{\Pc}_{\on{out}}$. But by \propref{coinv}, the space of
coinvariants of $\V^n_{k,c}|_{R}$ with respect to $\gtil_{\on{out}}$
is equal to the space of coinvariants of $\V^n_{k,c}|_{R}$ with
respect to $U_{X\bs x}^{\Pc}(\V_{k,c,x})$. This implies the statement
of the theorem, because according to the discussion of
\secref{reminder}, $U_{X\bs x}^{\Pc}(\V_{k,c,x})$ contains the image
of $\on{Vect}(X\bs x) \subset U_{X\bs x}(\Vir_{c_k,x})$ as a Lie
subalgebra.

\subsubsection{Descent to $\M_g$.}

The homomorphism of \thmref{descent to coinvariants} may be used to
describe differential operators on the moduli of unmarked curves and
bundles, $\Pi:\Bun_{G,g}\to\M_g$. To do so we first localize the
corresponding ``vacuum'' representations $Vir_{c_k}$ and $V_{k,c}$ on
$\Bun_{G,g,1}$ and $\Mgo$ and then descend.

\subsubsection{\bf Corollary.}\label{descent to curves}
{\em There is a canonical homomorphism $\D_{c_k}\longrightarrow
\Pi_*\D_{k,c}$ of sheaves of algebras on $\Mg$.}

\subsubsection{Proof.}
We will construct a morphism of sheaves $\Pi^*\D_{c_k}\to \D_{k,c}$ on
$\Bun_{G,g}$, which gives rise to the desired morphism on $\Mg$ by
adjunction. This map furthermore is constructed by descent from
$\Bun_{G,g,1}$.  By \propref{pullback of diffops} and
\propref{pullback of diffops for bundles}, the sheaves
$\Delta_{c_k}(Vir_{c_k})$ and $\Delta_{k,c}(V_{k,c})$ on $\Mgo$ and
$\Bun_{G,g,1}$ are identified with $\pi^*\D_{c_k}$ and
$\pi_G^*\D_{k,c}$. Moreover these identifications are {\em horizontal}
with respect to the relative connections inherited by the pullback
sheaves. We now need to construct a map between the pullback of
$\pi^*\D_{c_k}$ to $\Bun_{G,g,1}$ and $\pi_G^*\D_{k,c}$ which is flat
relative to $\pi_G$ and hence descends to $\Bun_{G,g}$.

This morphism $\Sug_{k,c}$ may be constructed directly following
\thmref{descent to coinvariants}, or by applying the change of vacuum
isomorphism \lemref{change of vacuum} to the homomorphism
$\Sug_{k,c}^1$. The connection along the curve $\pi_G:\Bun_{G,g,1}\to
\Bun_{G,g}$ is deduced (by passing to subquotients) from the action of
$L_{-1} \in \DerO$ on the vacuum modules $Vir_{c_k}$ and $V_{k,c}$,
which are compatible under the Segal-Sugawara homomorphism (see
\propref{Sugawara embedding}).  Hence the morphism is horizontal along
this curve and descends to $\Bun_{G,g,1}$.

The morphism $\Sug_{k,c}$ on $\Mg$ is a homomorphism, since it is
obtained by reduction (by change of vacuum, see \lemref{change of
vacuum}) from the corresponding homomorphism $\Sug^1_{k,c}$ on $\Mgo$.

\subsection{Tensor product decomposition.}

By applying localization to the
isomorphism of $(\gtil,\Gtil(\Oo))$--modules
$$\sug_{k,c}:V_k\ot\Vck \overset{\sim}\longrightarrow V_{k,c}$$ of
\propref{Sugawara embedding}, we obtain an isomorphism
$$\Delta\sug_{k,c}:\Delta(V_k\ot\Vck) \overset{\sim}\longrightarrow
\Delta(V_{k,c}).$$ However, the functors of coinvariants and
localization are not tensor functors. Since $V_k$ is a
$\ghat$--submodule of $V_{k,c}$, we obtain a natural map $\Delta(V_k)
\to \Delta(V_{k,c})$.  The localization $\Delta(V_k)$ is naturally
identified with the sheaf $\D_{k/\M}\subset \D_{k,c}$ of {\em
relative} differential operators.  We now use the homomorphism
$\Sug_{k,c}$ to lift the Virasoro operators to $\D_{k,c}$:

\subsubsection{\bf Theorem.}\label{quantum product}
{\em There is an isomorphism of sheaves on $\M_{g,1,2n}$,
$$\Pi_*^n\D_{k/\M} \und{\Oo}\ot \D_{c_k}\cong \Pi_*^n \D_{k,c},$$
compatible with the inclusions of the two factors as subalgebras of
$\Pi_*^n\D_{k,c}$.}

\subsubsection{Proof.}
Consider the isomorphism of $(\gtil,\Gtil_n(\Oo))$--modules
$$\sug_{k,c}^n:V_k^n \und{\C}\ot Vir_{c_k}^{2n} \longrightarrow
V_{k,c}^n$$ of \propref{Sugawara embedding}. Twisting by the
$\Gtil_n(\Oo)$--torsor $\wh{\Bun}_{G,g}\to \Bun_{G,g,1,n}$ we obtain
an isomorphism
$$\V_k^n \und{\Oo}\ot \Vir_{c_k}^{2n} \longrightarrow \V_{k,c}^n$$
of sheaves.  This morphism restricts to $\Sug_{k,c}^n$ on
$\vac\ot\Vir_{c_k}^{2n}$, and hence descends to an isomorphism
$$\Pi_*^n\V_k^n \und{\Oo}\ot \Vir_{c_k}^{2n}\longrightarrow \Pi_*^n
\V_{k,c}^n.$$ Thus we see that $\Pi_*^n \V_{k,c}^n$ is generated by
its two subalgebras $\Pi_*^n \V_k^n$ and $\Vir_{c_k}^n$.

The inclusion $V_k\to V_{k,c}$ of $\ghat$--modules gives rise to a
natural algebra inclusion
$$\Delta(V_k^n)=\D_{k/\M} \longrightarrow \Delta(V_{k,c}^n)=\D_{k,c}$$
from the sheaf of twisted differential operators {\em relative} to the
moduli of curves $\M_{g,1,2n}$ to the full sheaf of twisted
differential operators. We define the map $$\Pi_*^n\D_{k/\M}
\und{\Oo}\ot \D_{c_k}\to \Pi_*^n\D_{k,c}$$ as one generated by the
homomorphisms from $\Pi_*^n\D_{k/\M}$ and $\D_{c_k}$ to
$\Pi_*^n\D_{k,c}$.

This map is surjective because it comes from the composition
$$
\V^n_k \und{\Oo}\ot \Vir^{2n}_{c_k} \to \V^n_{k,c} \to
 \Delta(V^n_{k,c}),
$$ which is surjective because the first map is an isomorphism by
\propref{Sugawara embedding} and the second map is surjective by
definition. It remains to check that this map is injective. We have a
commutative diagram
\begin{equation*}
\begin{array}{ccc}
\V^n_k \und{\Oo}\otimes \Vir^{2n}_{c_k}
&\overset{\sim}\longrightarrow& \V^n_{k,c} \\ \downarrow&&\downarrow\\
\Delta(V^n_k) \und{\Oo}\otimes \Delta(Vir^{2n}_{c_k})
&\longrightarrow& \Delta(V^n_{k,c})
\end{array}
\end{equation*}
For an $R$--point of $\Bun_{G,g,1,n}$ as in the proof of
\thmref{descent to coinvariants}, the kernel of the right vertical map
over $\on{Spec} R$ is the image of the action of the Lie algebra
$\gtil^\Pc_{\on{out}}$ in $\V^n_{k,c}|_R \simeq \V_{k}|_R \ot
\Vir^{2n}_{c_k}|_R$. But
$$\gtil^\Pc_{\on{out}}\cdot (\V_{k}|_R \ot \Vir^{2n}_{c_k}|_R) \subset
(\g^\Pc_{\on{out}}\cdot \V_{k}|_R) \ot \Vir^{2n}_{c_k}|_R + \V_{k}|_R
\ot (\gtil^\Pc_{\on{out}}\cdot \Vir^{2n}_{c_k}|_R),$$ and so the
right hand side is in the kernel of the projection
$\V_{k}|_R \ot \Vir^{2n}_{c_k}|_R\to (\Delta(V^n_k) \ot
\Delta(\Vck))|_R$. This proves that our map is injective and hence an
isomorphism.

\subsubsection{Remark.} The theorem may be applied to
describe differential operators on the moduli of curves and bundles
without level structure, following \thmref{descent to
curves}. However, it is known (\cite{Hecke}) that for $k\neq -\hv$,
there are no non-constant global twisted differential operators on
$\Bun_G(X)$, in other words, $\Pi_*\D_{k/\M}=\Oo_{\M_g}$. Thus in this
case the pushforward of twisted differential operators $\D_{k,c}$ to
$\M_g$ is simply identified with the differential operators $\D_{c_k}$
on $\M_g$.

\subsection{Heat operators and projectively flat connections.}

Let $\pi:M\to S$ be a smooth projective morphism with connected fibers
and $\Ll\to M$ a line bundle. Let $\D_\Ll$ be the sheaf of twisted
differential operators on $\Ll$.

\subsubsection{Definition.} A heat operator on
$\Ll$ relative to $\pi$ is a lifting of the identity map
$\on{id}:\Theta_S \to \Theta_S$ to a sheaf homomorphism ${\mc H}:
\Theta_S\to \pi_*\D_{\Ll,\leq 1_S}$ to differential operators, which
are of order one along $S$, such that the corresponding map $\Theta_S
\to \pi_*\D_{\Ll,\leq 1_S}/\Oo_S$ is a Lie algebra homomorphism.

Suppose that the sheaf $\pi_* \Ll$ is locally free, i.e., is a sheaf
of sections of a vector bundle on $S$. Then a heat operator gives rise
to a projectively flat connection on this vector bundle (see
\cite{W,BK,Fa}).

Consider the morphism $\Pi: \Bun_{G,g} \to \M_g$. Recall that
$\Cc=\Ll_{k,0}$ is the line bundle on $\Bun_{G,g}$ whose restriction
to each $\Bun_{G}(X)$ is the ample generator of the Picard group of
$\Bun_{G}(X)$. For any $k \in \Z_+$ the sheaf $\Pi_* \Cc^k$ is locally
free, and its fiber at a curve $X\in\M_g$ is the vector space
$\on{H}^0(\Bun_G(X),\Cc^k)$ of {\em non-abelian theta functions} of
weight $k$. It is well-known that the corresponding vector bundle
(which we will also denote by $\Pi_* \Cc^k$), possesses a projectively
flat connection \cite{Hitchin connection,Fa,BK}. This connection may
be constructed using a heat operator on $\Cc^k$.

\subsubsection{\bf Theorem.} (\cite{BK})
{\em For any $k\geq 0$ the sheaf $\Pi_*\Cc^k$ on $\M_g$ possesses a
unique flat projective connection given by a heat operator on $\Cc^k$.}

\subsubsection{} The existence and uniqueness of this projective connection
are deduced for $\Ll$ satisfying the vanishing of the composition
$$\Theta_S\to R^1\pi_*\D_\Ll \to R^1\pi_*\D_{\Ll/S}$$ and the
identification $\pi_*\D_{\Ll/S}=\Oo_S$.  The connection was
also constructed by Hitchin \cite{Hitchin connection} and Faltings
\cite{Fa} by different means.

The connection on non-abelian theta functions has been explicitly
identified in \cite{Laszlo} with the connection on the (dual of) the
space of conformal blocks for the basic integrable representation
$L_k(\ghat)$ of $\ghat$ at level $k \in \Z_+$. The latter connection,
known as the KZB or WZW connection (and specializing to the
Knizhnik--Zamolodchikov connection in genus zero), is defined
following the general procedure of \cite{book}, Ch. 16. Namely, for
any conformal vertex algebra $V$ we defined a twisted $\D$--module of
coinvariants on $\M_g$. Its fibers are the spaces of coinvariants
$H(X,x,V)$ (see \secref{reminder}), and the action of the sheaf $\D_c$
of differential operators (where $c$ is the central charge of $V$)
comes from the action of the Virasoro algebra on $V$. In our case, we
take as $V$ the module $L_k(\ghat)$. For any positive integer $k$ this
module is a conformal vertex algebra with central charge $c(k) =
k\dim\g/(k+h^\vee)$, with the conformal structure defined by the
Segal-Sugawara vector. In this case the sheaf of coinvariants is
locally free as an $\Oo$--module, and so it is the sheaf of sections
of a vector bundle with a projectively flat connection. The sheaf of
sections of its dual vector bundle (whose fibers are the spaces of
conformal blocks of $L_k(\ghat)$) is therefore also a twisted
$\D$--module, more precisely, a $\D_{-c(k)}$--module. The
corresponding projectively flat connection on the bundle of conformal
blocks is the KZB connection.

On the other hand, in \thmref{descent to coinvariants} we produced
homomorphisms $\Sug_{k,0}^n:
\D_{-c(k)}\longrightarrow\Pi_*^n\D_{k,0}$. Their restrictions
$\Sug_{k,0}|_{\leq 1}$ to $(\D_{-c(k)})_{\leq 1}$ give us heat
operators on $\Cc^k=\Ll_{k,0}$. But the sheaf $\Pi_*\Cc^k$ is
isomorphic to the sheaf of conformal blocks on $\Mgo$ corresponding to
$L_k(\ghat)$ \cite{BL,KNR,Fa2,Tel}. Under this identification, the
projectively flat connection on $\Pi_*\Cc^k$ obtained from the heat
operators $\Sug_{k,0}|_{\leq 1}$ tautologically coincides with the KZB
connection, because both connections are constructed by applying the
Segal-Sugawara construction. Thus, we obtain

\subsubsection{\bf Proposition.}\label{heat operators}
{\em
\begin{enumerate}
\item For every $k\in\Z_+$, the heat operator defining the
projectively flat connection on the sheaf $\Pi_*\Cc^k$ of non-abelian
theta functions over $\Mg$ is given by the restriction of the
Segal-Sugawara map $\Sug_{k,0}|_{\leq 1} : (\D_{-c(k)})_{\leq 1}\to
\Pi_*\D_{k,0}$.

\item For any $n$, $\Sug_{k,0}^n|_{\leq 1}$ gives heat operators
defining a projectively flat connection on $\Pi^n_*\Ll_{k,c(k)}$ over
$\M_{g,1,2n}$.

\item
Under the identification between $\Pi_*\Cc^k$ and the sheaf of
conformal blocks on $\Mgo$ corresponding to the integrable
representation of $\ghat$ of level $k$ and highest weight $0$, the KZB
connection on the sheaf of conformal blocks is given by the heat
operators $\Sug_{k,0}|_{\leq 1}$ on $\Cc^k$.
\end{enumerate}
}

In the same way we obtain heat operators for any level structure
(where the morphisms $\Pi^n$ are no longer projective). Note also that
we can replace $\Cc^k$ by $\Ll_{k,c}$ for any integer $c$, which have
isomorphic restrictions to $Bun_G(X)$ for any $X$, and obtain
analogous projectively flat connections.

\section{Classical Limits.}\label{classical limits}

In this section we describe the limits of the (suitably rescaled)
Segal-Sugawara operators at the critical level $k=-\hv$, and in the
classical limit $k,c\to \infty$. In the former case the algebra of
twisted differential operators $\D_{c_k}$ degenerates into the Poisson
algebra of functions on the space of curves with projective structure,
and the construction gives twisted differential operators on
$Bun_{G,g}$ which are {\em vertical} (i.e. preserve $Bun_G(X)$ for
fixed $X$) and {\em commute}. We identify these operators with the
quadratic part of the Beilinson--Drinfeld quantization of the Hitchin
system. When the level and charge become infinite, both sides of the
construction become commutative (Poisson) algebras, and the
Segal--Sugawara construction is interpreted as a map from the moduli
of extended connections to the moduli of projective structures or
quadratic differentials. We interpret this map as defining a
symplectic connection over the moduli of pointed curves on the moduli
spaces of connections with arbitrary poles. We also sketch the
interpretation of this connection as a new Hamiltonian form of the
equations of isomonodromic deformation.

\subsection{The Critical level.}
Recall that we assume throughout that the group $G$ is
simply-connected.

By the general formalism of localization (\propref{twists are
algebras}), we have an algebra homomorphism $(V_k^n)^{G(\Oo)}\to
\Gamma(\Bun_{G}(X,x,n),\D_k)$ from the $G(\Oo)$--invariants of the
vacuum (which are the endomorphisms of the vacuum representation) to
global differential operators on the moduli of bundles.  For general
$k$ this gives only scalars, but at the critical level $k=-\hv$ the
space $V_k^{G(\Oo)}$ becomes very large.  In \cite{FF} (see also
\cite{Fr1,Fr2}), Feigin and Frenkel identify the algebra
$V_{-\hv}^{G(\Oo)}$ canonically with the ring of functions $\C[{\mc
Op}_{G^{\vee}}(D)]$ on the space of {\em opers} on the disc, for the
Langlands dual group $G^{\vee}$ of $G$.  Opers, introduced in
\cite{Hecke} by Beilinson and Drinfeld, are $G^{\vee}$--bundles
equipped with a Borel reduction and a connection, which satisfies a
strict form of Griffiths transversality.  For $G^{\vee}=PSL_2$ (so
$G=SL_2$), opers are identified with projective structures $\Proj(X)$.
In fact for arbitrary $G$ there is a natural projection
$\Op_{G^{\vee}}(X)\to \Proj(X)$, which identifies the space of opers
with the affine space for the vector space $\on{Hitch}_G(X)$ induced
from the affine space $\Proj(X)$ for quadratic differentials
$\on{H}^0(X,\Omega_X^{\ot 2})\subset \Hitch_G(X)$. (Similar
considerations apply to any level vacuum representations $V_k^n$,
replacing regular opers by opers with singularities, which form an
affine space for the meromorphic version of the Hitchin space.)

Thus there is a homomorphism
$$\C[\Op_{G^{\vee}}(D)]\longrightarrow \Gamma(\Bun_G(X),\D_{-\hv}).$$
Beilinson and Drinfeld show that this homomorphism factors through
functions on {\em global} opers $\C[\Op_{G^{\vee}}(X)]$, is
independent of the choice of $x\in X$ used in localization, and gives
rise to an {\em isomorphism}
$$\C[\Op_{G^{\vee}}(X)] \cong \Gamma(\Bun_G(X),\D_{-\hv}).$$
$$\Oo(\Op_{G^{\vee}}(\Xf))\to \Pi_*\D_{k/\M}.$$ We wish to compare the
restriction of this homomorphism to the subalgebra $\C[\Proj(X)]$ of
$\C[\Op_{G^{\vee}}(X)]$ with the critical level limit of the
Segal-Sugawara construction. Let $\Op_{G^{\vee}}\to \Mg$ denote the
moduli stack of curves with $G^{\vee}$--opers.

Recall the notation $\mu_\g=\hv \on{dim}\g$.


\subsubsection{\bf Theorem.}\label{critical theorem}
{\em
\begin{enumerate}
\item For $c\in\C$, the homomorphism $(k+\hv)\Sug_{k,c}^n$ is regular
at $k=-\hv$, defining an algebra homomorphism
  $$\ol{\Sug}_{-h^{\vee},c}^n:\Oo(T^*_{\mu_\g}\M_{g,1,2n})
\longrightarrow \Pi_*^n\D_{-h^{\vee},c}.$$
\item The homomorphism $\ol{\Sug}_{-\hv,c}$ is the restriction of the
Beilinson--Drinfeld homomorphism to projective structures: we have a
commutative diagram
\begin{equation*}
\begin{array}{ccc}
\Oo(\Proj^{\mu_\g}_g) & \stackrel{\ol{\Sug}_{-\hv,c}}{\longrightarrow} &
\Pi_*\D_{-\hv,c}\\ \downarrow && \uparrow \\ \Oo(\Op_{G^{\vee}}) &
\stackrel{\on{BD}}{\longrightarrow} & \Pi_*\D_{-\hv/\Mg}
\end{array}
\end{equation*}
\end{enumerate}
}

\subsubsection{Proof.}
As described in \secref{critical level} (see \propref{critical
Sugawara}), the rescaled Segal--Sugawara operators
$\ol{\Sug}_m=(k+\hv)\Sug_m$ are regular at $k=-\hv$, and define a
homomorphism of vertex algebras. We may now repeat the constructions
of \secref{Sugawara section} leading to \thmref{descent to
coinvariants} for the rescaled Sugawara operators. Note that the
classical vacuum representations $\ol{Vir}_{\mu_\g}$ are (commutative)
vertex algebras -- the vertex Poisson structure is not used in the
definition of coinvariants -- and for commutative vertex algebras, the
comparison of coinvariants for a generating set and the algebra it
generates is obvious.  Identifying the localization of
$\ol{Vir}_{\mu_\g}^{2n}$ with
$\Oo(T^*_{\mu_\g}\M_{g,1,2n})=\Oo(\Proj^{\mu_\g}_{g,1,2n})$, we obtain
the first assertion.

Moreover, by \secref{critical level}, the critical Segal-Sugawara
homomorphism factors as follows:
$$\ol{\Sug}_{-\hv,c}^n:\ol{Vir}_{\mu_\g}^{2n}\hookrightarrow
(V_{-\hv}^n)^{\Gn}\hookrightarrow (V_{-\hv,c}^n)^{\Gn}$$ (since
$\ol{\Sug}=-S_{-2}$ lands inside $V_{-h^{\vee}}\subset
V_{-h^{\vee},c}$).  It follows that the localized map on twisted
differential operators also factors through the localization of
$V_{-\hv}$, which is the sheaf of relative differential operators.
Thus the critical Segal-Sugawara construction is part of the
localization of the $G(\Oo)$--invariants in the vacuum representation,
giving the Beilinson--Drinfeld operators.

\subsection{Infinite level.}

Recall from \secref{moduli spaces} (see \propref{conn as TCB}) that
the $\Cc$--twisted cotangent bundle of the moduli stack of bundles
$\Bun_G(X)$ is isomorphic to the moduli stack $\Conn_G(X)$ of bundles
with regular connections on $X$, while the twisted cotangent bundle of
the moduli of bundles with $n$--jet of trivialization $\Bun_G(X,x,n)$
is the moduli $\Conn_G(X,x,n)$ of connections having at most $n$--th
order poles at $x$. Similarly, we may consider moduli
$\Bun_G(X,x_1,\dots,x_m,n_1,\dots,n_m)$ of bundles with several marked
points and jets of coordinates, whose twisted cotangent bundles are
identified with moduli of connections with poles of the corresponding
orders. As elsewhere, we restrict to the one--point case for
notational simplicity though all constructions carry over to the
multipoint case in a straightforward fashion.

By virtue of their identification as twisted cotangent bundles, the
moduli of meromorphic connections on a fixed curve carry canonical
(holomorphic) symplectic structures.  As we vary $X$ and $x$, the
moduli stack $\Conn_{G,g,1,n}$ forms a relative twisted cotangent
bundle to $\Bun_{G,g,1,n}$ over $\M_{g,1,2n}$. To obtain a symplectic
variety, we consider the {\em absolute} twisted cotangent bundles of
$\Bun_{G,g,1,n}$. The twisted cotangent bundle corresponding to the
line bundle $\Ll_{\la,\mu}$ is the moduli of {\em extended
connections}
$\Conn_{Bun_{G,g,1,n}}(\Ll_{\la,\mu})=\ExConn^{\la,\mu}_{G,g,1,n}$
defined in \secref{moduli of curves and bundles}. In the limit of
infinite level the sheaves of twisted differential operators on
$\Bun_{G,g,1,n}$ and $\M_{g,1,2n}$ degenerate to the commutative (and
hence Poisson) algebra of functions on the moduli of extended
connections and projective structures, respectively. Therefore it is
convenient to reinterpret the infinite limit of the Segal-Sugawara
homomorphism as a morphism between these moduli spaces (rather than a
homomorphism between the corresponding algebras of functions), and
examine its Poisson properties.

As in \secref{infinite limit}, we let $k,c\to \infty$ by introducing
auxiliary parameters $\la,\mu$ with
$\dfrac{\la}{k}=\dfrac{\mu}{c}$. Then the homomorphism
$\dfrac{\la^2}{k}\Sug_{k,c}^n$ is regular as $k,c\to\infty$, as is
(for $\la\neq 0$) the isomorphism $\sug^n_{k,c}:V_k^n\ot
Vir_{c_k}^{2n}\to V_{k,c}^n$. We thus obtain classical limits of
\thmref{descent to coinvariants} and \thmref{quantum product}.

Let $\on{Quad}_{g,1,n}=T^*\M_{g,1,n}$ and
$\Higgs_{G,g,1,n}=T^*_{/\M_{g,1,2n}}\Bun_{G,g,1,n}$ denote the moduli
stacks of curves equipped with a quadratic differential and bundles
with Higgs field, having at most $n$--th order pole at the marked
point, respectively. Let $\on{Hitch}_{G,g,1,n}$ denote the target of
the Hitchin map for $\Higgs_{G,g,1,n}$, i.e.,
$$
\on{Hitch}_{G,g,1,n} = \bigoplus_{i=1}^\ell \Gamma(X,\Omega(nx)^{d_i+1}),
$$
where $\ell = \on{rank} \g$ and $d_i$ is the $i$th exponent of $\g$.

Recall from \secref{infinite limit} that a commutative algebra
homomorphism of Poisson algebras (and the corresponding morphism of
spaces) is called $\la$--{\em Poisson} if it rescales the Poisson
bracket by $\la$.

\subsubsection{\bf Theorem.}\label{classical theorem}
{\em
\begin{enumerate}
\item
For every $\la,\mu\in\C$ there is a $\la$--Poisson homomorphism
$$\ol{\Sug}_{\la,\mu}^n:\Oo(T^*_{\la\mu}\M_{g,1,2n})\longrightarrow
\Pi_*^n\Oo(T^*_{\la,\mu}\Bun_{G,g,1,n})$$ Equivalently, for every
$\la,\mu\in\C$ there is a $\la$--Poisson map
  $$\Pi_{\la,\mu}^n:\ExConn^{\la,\mu}_{G,g,1,n}\to
\Proj^{\la\mu}_{g,1,2n}$$ lifting
$\Pi^n:\Bun_{G,g,1,n}\to\M_{g,1,2n}$.
\item For $\la\neq 0$, there is a canonical (non-affine) product
decomposition over the moduli of curves
$$\ExConn^{\la,\mu}_{G,g,1,n}\cong
\Conn_{G,g,1,n}^{\la} \und{\M_{g,1,2n}}\times \Proj^{\la\mu}_{g,1,2n}.$$
\item For $\la=\mu=0$, $\Pi_{\la,\mu}^n$ factors through the quadratic
Hitchin map
\begin{equation*}
\begin{array}{ccc}
T^*Bun_{G,g,1,n} & \stackrel{\Pi^n_{0,0}}{\longrightarrow} &
\on{Quad}_{g,1,2n}\\ \downarrow && \uparrow \\ \Higgs_{G,g,1,n} &
\stackrel{\on{Hitch}}{\longrightarrow} & \on{Hitch}_{G,g,1,n}
\end{array}
\end{equation*}
\end{enumerate}
}

\subsubsection{Proof.}
As in \thmref{critical theorem}, the construction of the classical
Segal-Sugawara homomorphism $\ol{\Sug}_{\la,\mu}$ is identical to the
proof of \thmref{descent to coinvariants}, appealing to
\propref{classical Sugawara} for the description of the rescaled
vertex algebra homomorphism. The localization of $\ol{Vir}_{\mu}^{2n}$
is $\Oo(T^*_{\mu}\M_{g,1,2n})$ and that of $\ol{V}_{\la,\mu}^n$ is
$\Oo(T^*_{\la,\mu}Bun_{G,g,1,2n})$, while the Poisson bracket is
rescaled by $\la$, by \propref{classical Sugawara}. The homomorphism
$\ol{\Sug}_{\la,\mu}^n$ of $\Oo_{\Bun_{G,g}}$--algebras defines, upon
taking $\on{Spec}$ over the moduli of curves, a morphism
$T^*_{\la,\mu}\Bun_{G,g,1,n}\to T^*_{\la\mu}\M_{g,1,2n}$ relative to
$\M_{g,1,2n}$, hence the geometric reformulation.

The morphism in (2) is given in components by the natural projection
\eqref{exconn to conn} from extended connections to connections and
the map $\Pi^n_{\la,\mu}$. The spaces $\ExConn^{\la,\mu}_{G,g,1,n}$
and $\Proj^{\la\mu}_{g,1,2n} \und{\M_{g,1,2n}}\times
\Conn^{\la}_{G,g,1,n}$ are both torsors for quadratic differentials
over $\Conn_{G,g,1,n}$. It follows from the explicit form of the
classical Segal-Sugawara vector that the resulting map of torsors over
the space of connections is an isomorphism (up to rescaling by
$\lambda$). Equivalently, the assertion (2) is the classical limit of
\thmref{quantum product} and may be proved identically, replacing
$\sug^n_{k,c}$ by its rescaled version.

For $\la=\mu=0$, the image of
$\ol{\Sug}^n_{0,0}:\ol{Vir}_0^{2n}\to\ol{V}_{0,0}^n$ lies strictly in
$\ol{V}_{0}^n\subset \ol{V}_{0,0}^n$. It follows that the image of
$\Oo(T^*\M_{g,1,2n})=\Oo(\on{Quad}_{g,1,2n})$ in
$\Oo(T^*Bun_{G,g,1,n})$ consists of functions which are pulled back
from $T^*_{/\M_{g,1,2n}}\Bun_{G,g,1,n}$. Therefore the map
$\Pi^n_{0,0}$ descends to
$T^*_{/\M_{g,1,2n}}\Bun_{G,g,1,n}$. Restricting to fibers over a fixed
curve $X$ we obtain a map $T^*\Bun_G(X,x,n)\to
T^*\M_{g,1,2n}|_X=\on{H}^0(X,\Omega^2(2nx))$.

The quadratic Hitchin map is the map $T^*\Bun_G(X,x,n)\to
\on{H}^0(X,\Omega^2(2nx))$ sending a Higgs bundle $(\Pc,\eta)$,
$\eta\in \on{H}^0(X,\g_{\Pc}\otimes\Omega(nx))$ to
$\dfrac{1}{2}\on{tr}\eta^2$.  Using the pairing on $L\g$, this is
identified with the map defined by the Segal-Sugawara vector
$$S_{-2}=\frac{1}{2}\sum_a J_{-1}^aJ_{a,-1}$$ as desired.

\subsubsection{Remark.} In \cite{thetas}, the product
decomposition (2) is described concretely in terms of kernel functions
along the diagonal (the necessary translation from vertex algebra
language to kernel language is described in Chapter 7 of \cite{book}).
Moreover, it is related, for $\mu\neq 0$, to the classical
constructions of projective structures on Riemann surfaces using theta
functions due to Klein and Wirtinger. Composing the projection
$\Pi_{\la,\mu}$ with the canonical meromorphic sections of the twisted
cotangent bundles over $\Bun_{G,g}$ given by non-abelian theta
functions (see also \cite{szego}), this gives a construction of
interesting rational maps from the moduli of bundles to the spaces of
projective structures, and more generally, opers.

Now we discuss various applications of the above results.

\subsubsection{The Bloch--Esnault--Beilinson connection.} \label{Bloch}
Let $\Xf\to S$ be a smooth family of curves, and $E$ a vector bundle
on $\Xf$. The vector bundle $E$ defines a line bundle $\Cc(E)$ on $S$,
the $c_2$--line bundle defined in \cite{Deligne} (whose Chern class is
the pushforward of $c_2(E)$) which differs from the determinant of
cohomology of $E$ by the $n$th power of the Hodge line bundle on
$S$. In fact, this line bundle is identified canonically with the
pullback of $\Cc$ to $S$ from the map $S\to Bun_{G,g}$ classifying $E$
and $X$ (see \cite{BK}). Bloch--Esnault and Beilinson \cite{BE}
construct a connection on $\Cc(E)$ from a regular connection
$\nabla_{\Xf/S}$ on $E$ relative to $S$. It is easy to see that this
connection can be recovered from \thmref{classical theorem} (2) (more
precisely, its straightforward version for $G=GL_n$ and unpointed
curves). Let $\la=1,\mu=0$, so that $\Ll_{\la,\mu}=\Cc$ and
$\Conn(\Cc)=\ExConn^{1,0}_{G,g}$ over $\Bun_{G,g}$. We have a
decomposition $\Conn(\Cc)=\Conn_{G,g} \und{\M_g}\times T^*\M_g$, and
hence a canonical lifting from $\Conn_{G,g}$ to $\Conn(\Cc)$, lying
over the zero section of $T^*\M_g$. When pulled back to $S$ by the
classifying map $S\to\Conn_{G,g}$ of $(E,\nabla_{\Xf/S})$, this gives
a connection on $\Cc(E)$ over $S$. The compatibility with the
construction of Bloch--Esnault--Beilinson follows for example from the
uniqueness of the splitting $\Conn_{G,g}\to\ExConn^{1,0}_{G,g}$, due
to the absence of one--forms on $\M_g$.

\subsubsection{The Segal-Sugawara symplectic connection.}
Suppose $f:N\to M$ is a smooth Poisson map of symplectic varieties.
It then follows that $N$ carries a flat symplectic connection over $M$,
i.e., a Lie algebra lifting of vector fields on $M$ to vector fields
on $N$ (defining a foliation on $N$ transversal to $f$), which
preserve the symplectic form on fibers (see \cite{Guil} for a
discussion of symplectic connections). Namely, we represent local
vector fields on $M$ by Hamiltonian functions, pull back to functions
on $N$ and take the symplectic gradient. The Poisson property of $f$
guarantees the flatness and symplectic properties of the
connection. (In the algebraic category, this defines the structure of
$\D$--scheme or crystal of schemes on $N$ over $M$, compatible with
the symplectic structure.)

The morphism $\Pi^n_{1,\mu}$ of \thmref{classical theorem} is such a
Poisson morphism of symplectic varieties, and hence defines a
symplectic connection (or crystal structure) on
$\ExConn^{1,\mu}_{G,g,1,n}$ over $\Proj^{\mu}_{g,1,2n}$. When $\mu=0$
we may reduce this connection as follows.  First we restrict to the
zero section $\M_{g,1,2n}\hookrightarrow \Proj^0_{g,1,2n}$ of the
cotangent bundle to obtain a connection on $\ExConn^{1,0}_{G,g,1,n}$
over $\M_{g,1,2n}$. Next, the connection respects the product
decomposition \thmref{classical theorem},(2) relative to $\M_{g,1,2n}$,
so that we obtain a flat symplectic connection on the moduli space
$\Conn_{G,g,1,n}$ of $G$--bundles with connections over the moduli
$\M_{g,1,2n}$ of decorated curves (with respect to the relative
symplectic structure):

\subsubsection{\bf Corollary.} {\em The projection $\Pi^n_{1,0}$ defines
a flat symplectic connection on the moduli stack $\Conn_{G,g,1,n}$
over $\M_{g,1,2n}$.}

\subsubsection{Time--dependent hamiltonians.}
The structure of flat symplectic connection on a relatively symplectic
$P\to M$ variety can also be encoded in a closed two--form $\Omega$ on
$P$, restricting to the symplectic form on the fibers. The connection
is then defined by the null--foliation of $\Omega$. If $P$ is locally
a product, then this structure can be encoded in the data of
Hamiltonian functions on $P$ that are allowed to depend on the
``times'' $M$. Thus, following \cite{M}, we may consider the data of
such a form as the general structure of time--dependent (or
non--autonomous) Hamiltonian system. In our case $P=\Conn_{G,g,1,n}$,
we thus have three equivalent formulations of a non-autonomous
Hamiltonian system, with times given by the moduli $\M_{g,1,2n}$ of
decorated curves: the Hamiltonian functions $\ol{\Sug}^n_{1,0}$ (or
more precisely the functions on $\Conn_{G,g,1,n}$ obtained from local
vector fields on $\M_{g,1,2n}$, considered as linear functions on
$T^*\M_{g,1,2n}$); the symplectic connection induced by $\Pi^n_{1,0}$; and
the two--form $\Omega$ on $\Conn_{G,g,1,n}$ obtained by restricting
the symplectic form on $\ExConn^{1,0}_{G,g,1,n}$ under the embedding
(as in \secref{Bloch}) along the zero section of $T^*\M_{g,1,2n}$.

\subsection{Isomonodromy.}
In this final section, we describe the algebraic definition of
isomonodromic deformation of arbitrary meromorphic connections over
the moduli of decorated curves and sketch its identification with the
Segal-Sugawara symplectic connection (as well as the compatibility
with the analytic iso--Stokes connections of \cite{Jimbo,Mal,Boalch}).
We thus obtain an algebraic time--dependent Hamiltonian description of
the isomonodromy equations. It also follows that the isomonodromy
Hamiltonians are classical limits of the heat operators defining the
KZB equations, and are non--autonomous deformations of the quadratic
Hitchin hamiltonians.

\subsubsection{Isomonodromic deformation.}
The moduli spaces of $G$--bundles with regular connection
$\Conn_{G,g}$ carry a flat connection (crystal structure) over the
moduli of curves, namely the connection of isomonodromic deformation
or nonabelian Gauss--Manin connection (see \cite{Si}). This connection
is a manifestation of the topological description of connections with
regular singularities on a Riemann surface as representations of the
fundamental group, which does not change under holomorphic
deformations. Namely, given a family $\Xf\to S$ of Riemann surfaces
and a bundle with (holomorphic, hence flat) connection on one fiber
$X$, there is a unique extension of the connection to nearby fibers so
that the monodromy representation does not change. The families of
connections relative to $S$ one obtains this way are uniquely
characterized as those families which admit an {\em absolute} flat
connection over $\Xf$, the total space of the family. This is captured
algebraically in the crystalline interpretation of flat connections: a
flat connection on a variety admits a unique flat extension to an
arbitrary nilpotent thickening of the variety. Thus algebraically one
defines families of flat connections relative to a base to be
isomonodromic if they may be extended to an absolute flat
connection. It follows that the moduli spaces of flat connections on
varieties carry a crystal structure over the deformation space of the
underlying variety -- the nonabelian Gauss-Manin connection of
Simpson \cite{Si}.

More generally, there is an algebraic connection (crystal or
$\D$--scheme structure) on the moduli stack of $\wh\Conn_{G,g,1,n}$ of
meromorphic connections over curves with formal coordinates, which is
the pullback of $\Conn_{G,g,1,n}$ to the moduli of pointed curves with
coordinates $\wh{\M}_{g,1}$. This connection is defined by ``fixing''
connections around their poles and deforming them isomonodromically on
the complement, combining the crystalline description of isomonodromy
and the Virasoro uniformization of the moduli of curves. Namely, given
a family $(\wt{X},\wt{x},\wt{t})\in \wh{\M}_{g,1}(S)$ of pointed
curves with formal coordinate over $S$, the spectrum of an Artinian
local ring, and a connection $(\Pc,\nabla)$ on the special fiber $X$
with pole at $x$ of order at most $n$, we must produce a canonical
extension of $(\Pc,\nabla)$ to $\wt{X}$. As explained in \cite{book},
Ch. 15, p.282, any such deformation $(\wt{X},\wt{x})$ of pointed
curves is given by ``regluing'' $X \bs x$ and the formal disc $D_x$
around $x$ using the action of $\AutK(S)$.

More precisely, we fix an identification of $\wt{X}\sm \wt{x}$ with
$(X\sm x)\times S$ and ``glue'' it to $D_x \times S$ using an
automorphism of the punctured disc over $S$. Now, given
$(\Pc,\nabla)$, we define an (absolute) flat connection on $\wt{X}$ by
extending our connection $(\Pc,\nabla)$ trivially (as a product, with
the trivial connection along the second factor) onto $(X \sm x)\times
S$ and onto $D_x \times S$. Note that since the connection on the
nilpotent thickening $\wt{X}\sm\wt{x} \simeq (X \sm x)\times S$ of
$X\sm x$ is flat, it is {\em uniquely} determined by $(\Pc,\nabla)$,
independently of the trivialization of the deformation -- is the
isomonodromic deformation of $(\Pc,\nabla)|_{X\sm x}$ over $S$ (in
particular the stabilizer of $(\wt{X},\wt{x})$ in $\AutK(S)$ does not
change the deformed connection.)  This defines the isomonodromy
connection on $\wh{\M}_{g,1}$.

\subsubsection{\bf Proposition.}\label{isomon}
{\em The Segal--Sugawara symplectic connection on $\Conn_{G,g,1,n}$
over $\M_{g,1,2n}$ pulls back to the isomonodromy connection on
$\wh\Conn_{G,g,1,n}$.}

\subsubsection{Proof.}
The compatibility between the two connections follows from the local
description of the isomonodromy connection through the action of
$\AutK$ on meromorphic connections, and the description of the
classical Segal-Sugawara operators as the corresponding Hamiltonians.
Namely, recall from \cite{book} that the ind--scheme
$\Conn_G(D^{\times})$ of connections on the trivial $G$--bundle on the
punctured disc is identified with $\ghat^*_1$, the level one
hyperplane in the dual to the affine Kac--Moody algebra. Moreover,
this identification is equivariant with respect to the natural actions
of $G(\K)$ and $\AutK$, describing the transformation of connections
under gauge transformations and changes of coordinates. These actions
are Hamiltonian, in an appropriately completed sense. Namely, the
Fourier coefficients of fields from the vertex Poisson algebra
$\ol{V}_{1}(\g)$ form the Lie subalgebra of {\em local functionals} on
$\ghat^*_1$ in the Poisson algebra of all functionals on
$\ghat^*_1$. By the Segal-Sugawara construction, this Lie algebra
contains as Lie subalgebra $\Der\K$ (as the Fourier coefficients
$\ol{L_n}$ of the classical limit of the conformal vector), thereby
providing Hamiltonians for the action of $\AutK$. The construction of
the isomonodromy connection on $\wh{\M}_{g,1}$ above is expressed in
local coordinates (as a flow on meromorphic connections on the disc at
$x$) by this action of $\AutK$ on connections. Likewise, the
Segal-Sugawara connection is defined by Hamiltonian functions, which
are reductions of the classical Segal-Sugawara operators to
$\Conn_{G,g,1,n}$. Thus the compatibility of the two connections is a
consequence of the local statement.

\subsubsection{Remark: analytic isomonodromy and Stokes data.}
In the complex analytic setting, one can extend to irregular
connections the description of regular--singular connections by
topological monodromy data, by introducing Stokes data describing the
transitions between asymptotic fundamental solutions to the connection
in different sectors. One can thereby define an ``iso--Stokes''
generalization of the isomonodromy equations, i.e., an analytic
symplectic connection on moduli of irregular holomorphic connections
(see \cite{Jimbo,Mal} and recently \cite{Boalch,Krich}). Thus, one has
an isomonodromic deformation with more ``times'' than just the moduli
of pointed curves (one may also vary the most singular term of the
connection).

It is easy to see our algebraic isomonodromic connection (with respect
to motions of pointed curves) agrees with the analytic one. Since we
can choose the gluing transformations in $\AutK$ used to describe a
given deformation to be {\em convergent} rather than formal, it
follows that we are not altering the isomorphism type of the
connection on a small analytic disc around $x$, and hence all Stokes
data are automatically preserved. Thus it follows that the iso--Stokes
connection on $\Conn_{G,g,1,n}$ is in fact algebraic.

\subsubsection{Conclusion: KZB, isomonodromy and Hitchin.}
It follows from \propref{isomon} and \thmref{classical theorem} that
the moduli spaces of meromorphic connections carry algebraic
isomonodromy equations, which form a time--dependent Hamiltonian
system. Moreover, the isomonodromy hamiltonians are non--autonomous
deformations of the quadratic Hitchin hamiltonians (see also
\cite{Krich}). Most significantly, it follows that the heat operators
$\Sug_{k,c}^n$ defining the KZB connection on the bundles of conformal
block quantize the isomonodromy Hamiltonians $\ol{\Sug}^n_{1,0}$. This
makes precise the picture developed in \cite{Ivanov} and \cite{LO} of
the non--stationary Schr\"odinger equations defining the KZB
connection on conformal blocks as quantizations of the isomonodromy
equations or non--autonomous Hitchin systems.  In genus zero we thus
generalize the result of \cite{Res,Har} that the KZ connection on
spaces of conformal blocks on the $n$--punctured sphere, viewed as a
system of multi--time--dependent Schr\"odinger equations, quantizes
the Schlesinger equation, describing isomonodromic deformation of
connections with regular singularities on the sphere, which itself is
a time--dependent deformation of the Gaudin system (the Hitchin system
corresponding to the $n$--punctured sphere). In genus one we obtain a
similar relation between KZB equations, the elliptic form of the
Painlev\'e equations and the elliptic Calogero--Moser system as in
\cite{LO}.

\end{document}